\documentclass[12pt]{article}%
\usepackage{amsfonts}
\usepackage{sw20bams}
\usepackage{amsmath}
\usepackage{amssymb}
\usepackage{graphicx}%
\providecommand{\U}[1]{\protect\rule{.1in}{.1in}}
\begin{document}

\title{$n^{\text{th}}$ Roots of $n^{\text{th}}$ Powers}
\author{Steven Finch}
\date{April 15, 2026}
\maketitle

\begin{abstract}
Seeking simple, efficient solutions of a matrix equation leads (quite
circuitously) to optimizing unimodular zerofree square matrices.
\ Canonicalizing such matrices under signed-permutation double action offers
an ideal application of GPUs (graphics processing units).

\end{abstract}

\footnotetext{Copyright \copyright \ 2026 by Steven R. Finch. All rights
reserved.}

Our story unfolds in 1879, when three (real)\ $2\times2$ matrix solutions of
$X^{3}=C^{3}$ were unveiled \cite{John-mtrx}, assuming%
\[%
\begin{array}
[c]{ccc}%
C=\left(
\begin{array}
[c]{cc}%
1 & -2\\
2 & -1
\end{array}
\right)  , &  & C^{3}=\left(
\begin{array}
[c]{cc}%
-3 & 6\\
-6 & 3
\end{array}
\right)  .
\end{array}
\]
Seeing this result, we wondered right away about larger exponents. \ Let $i$
denote the imaginary unit. \ For odd integer $n\geq3$, let $\omega=\exp
(2\pi\,i/n)$, $\xi=\exp(\pi\,i/3)$ and $\eta=\exp(\pi\,i/(2n))$. \ Define
matrices%
\[
X_{n,j,k}=\frac{(-1)^{(n+1)/2}}{i}\left(
\begin{array}
[c]{cc}%
-\dfrac{1}{\xi}\,\eta\,\omega^{j}+\xi\,\dfrac{1}{\eta}\,\omega^{k} &
\eta\,\omega^{j}-\dfrac{1}{\eta}\,\omega^{k}\medskip\\
-\eta\,\omega^{j}+\dfrac{1}{\eta}\,\omega^{k} & \xi\,\eta\,\omega^{j}%
-\dfrac{1}{\xi\,\eta}\,\omega^{k}%
\end{array}
\right)
\]
for $0\leq j,k\leq n-1$. \ There exist nine (complex)\ solutions of
$X^{3}=C^{3}$, three of which are indeed real:%
\[%
\begin{array}
[c]{ccccc}%
X_{3,0,0}=\left(
\begin{array}
[c]{cc}%
1 & 1\\
-1 & 2
\end{array}
\right)  , &  & X_{3,1,2}=\left(
\begin{array}
[c]{cc}%
-2 & 1\\
-1 & -1
\end{array}
\right)  , &  & X_{3,2,1}=\left(
\begin{array}
[c]{cc}%
1 & -2\\
2 & -1
\end{array}
\right)  .
\end{array}
\]
A typical solution possessing nonzero imaginary part is%
\[
X_{3,1,0}=\left(
\begin{array}
[c]{cc}%
-\dfrac{1}{2}-i\dfrac{\sqrt{3}}{2} & 1+i\sqrt{3}\\
-1-i\sqrt{3} & \dfrac{1}{2}+i\dfrac{\sqrt{3}}{2}%
\end{array}
\right)  =\omega^{2}X_{3,2,1}.
\]
More generally, there exist $n^{2}$ solutions of $X^{n}=C^{n}$. \ For $n=5$,
there are five real solutions, including%
\[%
\begin{array}
[c]{cc}%
X_{5,3,2}=\left(
\begin{array}
[c]{cc}%
\dfrac{1}{4}\left(  -1-\sqrt{5}+\sqrt{30-6\sqrt{5}}\right)  & \dfrac{1}%
{2}\left(  1+\sqrt{5}\right) \\
\dfrac{1}{2}\left(  1-\sqrt{5}\right)  & \dfrac{1}{4}\left(  1+\sqrt{5}%
+\sqrt{30-6\sqrt{5}}\right)
\end{array}
\right)  , & X_{5,1,4}=C.
\end{array}
\]
Algebraic expressions regrettably become too complicated beyond this point.

The situation for even integer $n\geq2$ is completely different:\ $C^{n}%
=(-3)^{n/2}I$, where $I$ is the $2\times2$ identity matrix \cite{Chdy-mtrx}.
\ Let $\omega$ be as before and%
\[
\kappa=\sqrt{3}\exp\left[  \left(  1-(-1)^{n/2}\right)  \frac{\pi\,i}%
{2n}\right]  .
\]
It turns out that there are infinitely many solutions of $Y^{n}=C^{n}$:%
\[
Y_{n,u,v}=\dfrac{1}{2}\left(
\begin{array}
[c]{cc}%
\kappa\left(  1+\omega\right)  +2u & v\left[  \kappa\left(  1-\omega\right)
+2u\right] \\
\dfrac{1}{v}\left[  \kappa\left(  1-\omega\right)  -2u\right]  & \kappa\left(
1+\omega\right)  -2u
\end{array}
\right)
\]
where $u$, $v$ are arbitrary parameters with $v\neq0$. \ For example,%
\[%
\begin{array}
[c]{ccc}%
Y_{2,-i\sqrt{3},2}=i\sqrt{3}\left(
\begin{array}
[c]{cc}%
-1 & 0\\
1 & 1
\end{array}
\right)  , &  & Y_{2,0,-i\sqrt{3}}=\left(
\begin{array}
[c]{cc}%
0 & -3\\
1 & 0
\end{array}
\right)  .
\end{array}
\]
The striking difference in behavior (odd $n$ versus even $n$) can be easily
explained. \ The matrix $C$ possesses distinct eigenvalues $\pm i\sqrt{3}$
(which is good) whereas the eigenvalues for $C^{2}$ coincide (which is bad).
\ We are interested in deriving more examples of the former (finite) type that
are relatively simple. \ It seems that the proper choice of integer
eigenvectors would help, as well as ensuring that integer eigenvalues do not
share the same absolute value. \ Our approach is experimental:\ we make no
claim of originality (as surely the literature contains examples that have
defied our detection). \ Insight and corrections from readers are always welcome.

\section{Integer $2\times2$ Matrices}

The narrative skips ahead 125 years, when sixteen solutions of $X^{4}=B^{4}$
were listed \cite{Chdy-mtrx}, assuming%
\[%
\begin{array}
[c]{ccc}%
B=\left(
\begin{array}
[c]{cc}%
-1 & 6\\
-2 & 6
\end{array}
\right)  , &  & B^{4}=\left(
\begin{array}
[c]{cc}%
-179 & 390\\
-130 & 276
\end{array}
\right)  .
\end{array}
\]
For integer $n\geq2$, define matrices%
\[%
\begin{array}
[c]{ccc}%
X_{n,j,k}=\left(
\begin{array}
[c]{cc}%
8\omega^{j}-9\omega^{k} & -12\omega^{j}+18\omega^{k}\\
4\omega^{j}-6\omega^{k} & -6\omega^{j}+12\omega^{k}%
\end{array}
\right)  , &  & 0\leq j,k\leq n-1
\end{array}
\]
where $\omega=\exp(2\pi\,i/n)$. \ Everyone can agree that $X_{n,j,k}$ here is
simpler than the preceding $X_{n,j,k}$. \ Four of the solutions are real; we
give two:%
\[%
\begin{array}
[c]{ccc}%
X_{4,0,0}=\left(
\begin{array}
[c]{cc}%
-1 & 6\\
-2 & 6
\end{array}
\right)  , &  & X_{4,2,0}=\left(
\begin{array}
[c]{cc}%
-17 & 30\\
-10 & 18
\end{array}
\right)
\end{array}
\]
as well as two solutions with nonzero imaginary part:%
\[%
\begin{array}
[c]{ccc}%
X_{4,1,2}=\left(
\begin{array}
[c]{cc}%
9+8i & -18-12i\\
6+4i & -12-6i
\end{array}
\right)  , &  & X_{4,3,2}=\left(
\begin{array}
[c]{cc}%
9-8i & -18+12i\\
6-4i & -12+6i
\end{array}
\right)  .
\end{array}
\]
The remaining solutions are the same as these, up to a root-of-unity factor
$\{i,-1,-i\}$.

From where do the expressions for $X_{n,j,k}$ arise?\ \ Diagonalization
provides%
\[
B=M\,\Lambda\,M^{-1}=\left(
\begin{array}
[c]{cc}%
2 & 3\\
1 & 2
\end{array}
\right)  \left(
\begin{array}
[c]{cc}%
2 & 0\\
0 & 3
\end{array}
\right)  \left(
\begin{array}
[c]{cc}%
2 & -3\\
-1 & 2
\end{array}
\right)
\]
thus%
\[
X_{n,j,k}=\left(
\begin{array}
[c]{cc}%
2 & 3\\
1 & 2
\end{array}
\right)  \left(
\begin{array}
[c]{cc}%
2\omega^{j} & 0\\
0 & 3\omega^{k}%
\end{array}
\right)  \left(
\begin{array}
[c]{cc}%
2 & -3\\
-1 & 2
\end{array}
\right)  ,
\]
The columns in the leftmost matrix $M$ are the eigenvectors $(2,1)$ and
$(3,2)$ of $B$. \ The determinant of $M$ is $1$; if we were to reverse the
order of eigenvalues, then $M$ would have determinant $-1$. \ 

A square integer matrix with determinant $\pm1$ is called \textbf{unimodular}.
\ Such matrices play a role in finding marginally simpler yet more
\textquotedblleft efficient\textquotedblright\ examples $A$. \ By this, we
mean reducing the size of $\left\Vert A\right\Vert $, where $\left\Vert
A\right\Vert $ denotes the absolute maximum over all entries in $A$. \ For
instance,%
\[
\left(
\begin{array}
[c]{cc}%
1 & 1\\
1 & 2
\end{array}
\right)  \left(
\begin{array}
[c]{cc}%
2 & 0\\
0 & 3
\end{array}
\right)  \left(
\begin{array}
[c]{cc}%
2 & -1\\
-1 & 1
\end{array}
\right)  =\left(
\begin{array}
[c]{cc}%
1 & 1\\
-2 & 4
\end{array}
\right)
\]
and the maximum $4$ improves upon the earlier maximum $6$ (for $B$).
\ Interestingly, if we change the eigenvalues $\{2,3\}$ to $\{1,2\}$, then%
\[
\left(
\begin{array}
[c]{cc}%
1 & 1\\
1 & 2
\end{array}
\right)  \left(
\begin{array}
[c]{cc}%
1 & 0\\
0 & 2
\end{array}
\right)  \left(
\begin{array}
[c]{cc}%
2 & -1\\
-1 & 1
\end{array}
\right)  =\left(
\begin{array}
[c]{cc}%
0 & 1\\
-2 & 3
\end{array}
\right)
\]
has maximum $3$ -- a further improvement -- but at the price of permitting a
zero entry in $A$. \ If we prohibit zeroes in $A$, then%
\[
\left(
\begin{array}
[c]{cc}%
1 & 2\\
1 & 3
\end{array}
\right)  \left(
\begin{array}
[c]{cc}%
1 & 0\\
0 & 2
\end{array}
\right)  \left(
\begin{array}
[c]{cc}%
3 & -2\\
-1 & 1
\end{array}
\right)  =\left(
\begin{array}
[c]{cc}%
-1 & 2\\
-3 & 4
\end{array}
\right)
\]
and the maximum returns to $4$. \ 

Let's stipulate that $\Lambda=\operatorname*{diag}(\lambda,\mu)$ where
$\lambda<\mu$ are positive integers. \ An invertible matrix $Z$ is called
\textbf{zerofree} if none of the entries in $Z$ and none of the entries in
$Z^{-1}$ are zero. \ (Involving both matrices in our definition is somewhat
unorthodox.) \ Several optimization problems come to mind:

\begin{enumerate}
\item[(i)] Find a $2\times2$ unimodular zerofree matrix $M$ such that
$\left\Vert \left(  M\;M^{-1}\right)  \right\Vert $ is minimal

\item[(ii)] Find $M$ as in (i) such that, given $\Lambda$, $\left\Vert
M\,\Lambda\,M^{-1}\right\Vert $ is minimal

\item[(iii)] Same as (ii) except $M\,\Lambda\,M^{-1}$ is additionally
constrained to be zerofree
\end{enumerate}

\noindent where $\left(  M\;M^{-1}\right)  $ denotes $2\times4$ concatenation.
\ We are unaware of any guiding theory here and have no choice but to apply a
brute-force algorithm. \ The solution of (i)\ is unsurprisingly $\tbinom
{1\;\;1}{1\;\;2}$; the solution of (iii) when $\Lambda=\operatorname*{diag}%
(1,2)$ is $\tbinom{1\;\;2}{1\;\;3}$. \ Investigating how (ii)\ and (iii)
depend on $\{\lambda,\mu\}$ is left to the reader.

Circling back to the main topic, given%
\[%
\begin{array}
[c]{ccc}%
A=\left(
\begin{array}
[c]{cc}%
1 & 1\\
-2 & 4
\end{array}
\right)  , &  & \text{with eigenvalues }\{2,3\}\text{,}%
\end{array}
\]
the solution of $X^{n}=A^{n}$ is%
\[%
\begin{array}
[c]{ccc}%
X_{n,j,k}=\left(
\begin{array}
[c]{cc}%
4\omega^{j}-3\omega^{k} & -2\omega^{j}+3\omega^{k}\\
4\omega^{j}-6\omega^{k} & -2\omega^{j}+6\omega^{k}%
\end{array}
\right)  , &  & 0\leq j,k\leq n-1.
\end{array}
\]
This is the simplest $n^{\text{th}}$ root (of an $n^{\text{th}}$ power)
expression we have yet seen. \ We must now confess that our phrase
\textquotedblleft the solution\textquotedblright\ is not quite
accurate:\ there are, in fact, $32$ solutions of (ii)\ for $M\,\Lambda
\,M^{-1}$, given $\Lambda=\operatorname*{diag}(2,3)$. \ It can be shown,
however, that only $16$ solutions are equivalent to $\tbinom{\,\,1\;\;\,1}%
{-2\;\;4}$, in the sense that one is obtained from the other via
pre-multiplication and post-multiplication by signed permutation matrices.
\ The remaining $16$ solutions are equivalent to%
\[%
\begin{array}
[c]{ccc}%
\tilde{A}=\left(
\begin{array}
[c]{cc}%
1 & 2\\
-1 & 4
\end{array}
\right)  , &  & \text{again with eigenvalues }\{2,3\}\text{;}%
\end{array}
\]
hence%
\[%
\begin{array}
[c]{ccc}%
\tilde{X}_{n,j,k}=\left(
\begin{array}
[c]{cc}%
4\omega^{j}-3\omega^{k} & -4\omega^{j}+6\omega^{k}\\
2\omega^{j}-3\omega^{k} & -2\omega^{j}+6\omega^{k}%
\end{array}
\right)  , &  & 0\leq j,k\leq n-1.
\end{array}
\]
The distinction between $X$ and $\tilde{X}$ is almost not worth mentioning.
\ But this little diversion highlights the advantage of solving (i)\ rather
than (ii)\ or (iii): the matrices $M=\tbinom{1\;\;1}{1\;\;2}$ \& $\tilde
{M}=\tbinom{-2\;\;1}{-1\;\;1}$, although they critically underpin inequivalent
$A$ \& $\tilde{A}$, are (in themselves) equivalent:%
\[
P\,M\,Q=\left(
\begin{array}
[c]{cc}%
0 & -1\\
-1 & 0
\end{array}
\right)  \left(
\begin{array}
[c]{cc}%
1 & 1\\
1 & 2
\end{array}
\right)  \left(
\begin{array}
[c]{cc}%
0 & -1\\
1 & 0
\end{array}
\right)  =\left(
\begin{array}
[c]{cc}%
-2 & 1\\
-1 & 1
\end{array}
\right)  =\tilde{M}.
\]
The signed permutation matrices $P$ \&\ $Q$ can alternatively be replaced by
$-P$ \&\ $-Q$, of course. \ We shall shift focus to (i) henceforth.

\section{Integer $3\times3$ Matrices}

From%
\[%
\begin{array}
[c]{ccc}%
Z=\left(
\begin{array}
[c]{ccc}%
1 & 0 & 0\\
1 & 1 & 0\\
1 & 1 & 1
\end{array}
\right)  \left(
\begin{array}
[c]{ccc}%
1 & 1 & 1\\
0 & 1 & 1\\
0 & 0 & 1
\end{array}
\right)  =\left(
\begin{array}
[c]{ccc}%
1 & 1 & 1\\
1 & 2 & 2\\
1 & 2 & 3
\end{array}
\right)  , &  & Z^{-1}=\left(
\begin{array}
[c]{ccc}%
2 & -1 & 0\\
-1 & 2 & -1\\
0 & -1 & 1
\end{array}
\right)
\end{array}
\]
we see that, while $Z$ is unimodular, it fails to be zerofree. \ It would
otherwise have been a nice analog of $\tbinom{1\;\;1}{1\;\;2}=\tbinom
{1\;\;0}{1\;\;1}\tbinom{1\;\;1}{0\;\;1}$.

There are $576$ solutions of (i), all of which are equivalent to $M$, where%
\[%
\begin{array}
[c]{ccc}%
M=\left(
\begin{array}
[c]{ccc}%
1 & 2 & 2\\
2 & 1 & 2\\
2 & 2 & 3
\end{array}
\right)  , &  & M^{-1}=\left(
\begin{array}
[c]{ccc}%
1 & 2 & -2\\
2 & 1 & -2\\
-2 & -2 & 3
\end{array}
\right)  .
\end{array}
\]
Let's stipulate that $\Lambda=\operatorname*{diag}(\lambda,\mu,\nu)$ where
$\lambda<\mu<\nu$ are positive integers. For eigenvalues $\{1,2,3\}$,\
\[
A=M\,\Lambda\,M^{-1}=\left(
\begin{array}
[c]{ccc}%
-3 & -6 & 8\\
-6 & -6 & 10\\
-8 & -10 & 15
\end{array}
\right)
\]
thus the solution of $X^{n}=A^{n}$ is
\[
X_{n,j,k,\ell}=\left(
\begin{array}
[c]{ccc}%
\omega^{j}+8\omega^{k}-12\omega^{\ell} & 2\omega^{j}+4\omega^{k}%
-12\omega^{\ell} & -2\omega^{j}-8\omega^{k}+18\omega^{\ell}\\
2\omega^{j}+4\omega^{k}-12\omega^{\ell} & 4\omega^{j}+2\omega^{k}%
-12\omega^{\ell} & -4\omega^{j}-4\omega^{k}+18\omega^{\ell}\\
2\omega^{j}+8\omega^{k}-18\omega^{\ell} & 4\omega^{j}+4\omega^{k}%
-18\omega^{\ell} & -4\omega^{j}-8\omega^{k}+27\omega^{\ell}%
\end{array}
\right)
\]
for $0\leq j,k,\ell\leq n-1$, where $\omega=\exp(2\pi\,i/n)$. \ As an example,
the four real $4^{\text{th}}$ roots of $A^{4}$, up to sign, are%

\[%
\begin{array}
[c]{ccccc}%
X_{3,0,0,0}=A, &  & X_{3,0,0,2}=\left(
\begin{array}
[c]{ccc}%
21 & 18 & -28\\
18 & 18 & -26\\
28 & 26 & -39
\end{array}
\right)  , &  & X_{3,0,2,0}=\left(
\begin{array}
[c]{ccc}%
-19 & -14 & 24\\
-14 & -10 & 18\\
-24 & -18 & 31
\end{array}
\right)  ,\medskip\\
&  & X_{3,0,2,2}=\left(
\begin{array}
[c]{ccc}%
5 & 10 & -12\\
10 & 14 & -18\\
12 & 18 & -23
\end{array}
\right)  . &  &
\end{array}
\]
Two papers we nearly overlooked \cite{Hau1-mtrx, Hau2-mtrx} contain a relevant
quote --\ \textquotedblleft Mathematical principles must not be hidden by
lengthy calculations\textquotedblright\ -- and examine topics like solving
linear systems, finding eigenvalues and Jordan normal forms with an eye on
pedagogical simplicity. \ 

\section{Integer $4\times4$ Matrices}

There are myriad solutions of (i), one-seventh of which are equivalent to
$M_{1}$, where%
\[%
\begin{array}
[c]{ccc}%
M_{1}=\left(
\begin{array}
[c]{cccc}%
1 & 1 & 1 & 2\\
1 & 1 & 2 & 1\\
1 & 2 & 2 & 2\\
2 & 1 & 2 & 2
\end{array}
\right)  , &  & M_{1}^{-1}=\left(
\begin{array}
[c]{cccc}%
-2 & -2 & 1 & 2\\
-2 & -2 & 2 & 1\\
1 & 2 & -1 & -1\\
2 & 1 & -1 & -1
\end{array}
\right)  ;
\end{array}
\]
two-sevenths of which are equivalent to $N_{2}$, where
\[%
\begin{array}
[c]{ccc}%
N_{2}=\left(
\begin{array}
[c]{cccc}%
1 & 1 & 1 & 2\\
1 & 2 & 2 & 2\\
1 & -1 & -2 & 1\\
2 & -1 & -2 & 2
\end{array}
\right)  , &  & N_{2}^{-1}=\left(
\begin{array}
[c]{cccc}%
-2 & 1 & -2 & 2\\
-2 & 2 & 2 & -1\\
1 & -1 & -2 & 1\\
2 & -1 & 1 & -1
\end{array}
\right)  ;
\end{array}
\]
and four-sevenths of which are equivalent to $N_{3}$, where
\[%
\begin{array}
[c]{ccc}%
N_{3}=\left(
\begin{array}
[c]{cccc}%
1 & 1 & 1 & 1\\
1 & 1 & 2 & 2\\
1 & 2 & 1 & 2\\
1 & 2 & -1 & 1
\end{array}
\right)  , &  & N_{3}^{-1}=\left(
\begin{array}
[c]{cccc}%
1 & 1 & -2 & 1\\
1 & -2 & 2 & -1\\
1 & -1 & 1 & -1\\
-2 & 2 & -1 & 1
\end{array}
\right)  .
\end{array}
\]
More precisely, choosing a matrix $X$ uniformly at random from the pool of
$129024$ minimal unimodular zerofree $4\times4$ matrices and defining
$\Xi=\operatorname*{CanonicalizeMatrix}[X]$, we obtain exact probabilities%
\[%
\begin{array}
[c]{ccccc}%
\mathbb{P}(\Xi=M_{1})=\dfrac{1}{7}, &  & \mathbb{P}(\Xi=N_{2})=\dfrac{2}{7}, &
& \mathbb{P}(\Xi=N_{3})=\dfrac{4}{7}.
\end{array}
\]
It is fascinating that $\left\Vert \left(  M\;M^{-1}\right)  \right\Vert =2$
here, but the norm of $3\times3$ unimodular zerofree matrices was $\geq3$ in
the preceding section. \ Why should increasing the dimension from $3$ to $4$
actually improve the efficiency of the best-case scenario? \ 

Things are about to become more intricate. \ From this point onward, we adopt
a convention:\ $M$ will be used for matrices with positive entries only,
whereas $N$ will be used for matrices containing at least one negative entry.

\section{Integer $5\times5$ Matrices}

Enumerating the minimal unimodular zerofree $5\times5$ matrices is unsolved,
but we compute that each matrix is equivalent to either
\[%
\begin{array}
[c]{ccc}%
M_{1}=\left(
\begin{array}
[c]{ccccc}%
1 & 1 & 1 & 2 & 2\\
1 & 1 & 2 & 1 & 2\\
1 & 2 & 2 & 2 & 2\\
2 & 1 & 2 & 2 & 2\\
2 & 2 & 2 & 2 & 1
\end{array}
\right)  , &  & M_{1}^{-1}=\left(
\begin{array}
[c]{ccccc}%
2 & 2 & -3 & -2 & 2\\
2 & 2 & -2 & -3 & 2\\
-3 & -2 & 3 & 3 & -2\\
-2 & -3 & 3 & 3 & -2\\
2 & 2 & -2 & -2 & 1
\end{array}
\right)
\end{array}
\]
or
\[%
\begin{array}
[c]{ccc}%
N_{2}=\left(
\begin{array}
[c]{ccccc}%
1 & 1 & 1 & 1 & 1\\
1 & 1 & 1 & 2 & -1\\
1 & 1 & 2 & 1 & -1\\
1 & 2 & 2 & 2 & -1\\
2 & 1 & 2 & 2 & -1
\end{array}
\right)  , &  & N_{2}^{-1}=\left(
\begin{array}
[c]{ccccc}%
-1 & -3 & -3 & 2 & 3\\
-1 & -3 & -3 & 3 & 2\\
1 & 2 & 3 & -2 & -2\\
1 & 3 & 2 & -2 & -2\\
1 & 1 & 1 & -1 & -1
\end{array}
\right)  ;
\end{array}
\]
demonstrating that $\left\Vert \left(  M\;M^{-1}\right)  \right\Vert =3$
although $\left\Vert M\right\Vert =2<\left\Vert M^{-1}\right\Vert $. \ We did
not see this \textquotedblleft divergent\textquotedblright\ behavior for
minimal $3\times3$ or $4\times4$ matrices.

\section{Integer $6\times6$ Matrices}

Enumerating the minimal unimodular zerofree $6\times6$ matrices is again
unsolved, but we compute that each matrix is equivalent to
\[%
\begin{array}
[c]{ccc}%
M_{1}=\left(
\begin{array}
[c]{cccccc}%
1 & 1 & 1 & 1 & 1 & 2\\
1 & 1 & 1 & 1 & 2 & 1\\
1 & 1 & 1 & 2 & 1 & 1\\
1 & 1 & 2 & 1 & 2 & 2\\
1 & 2 & 1 & 2 & 1 & 2\\
2 & 1 & 1 & 2 & 2 & 1
\end{array}
\right)  , &  & M_{1}^{-1}=\left(
\begin{array}
[c]{cccccc}%
-1 & -2 & -2 & 1 & 1 & 2\\
-2 & -1 & -2 & 1 & 2 & 1\\
-2 & -2 & -1 & 2 & 1 & 1\\
1 & 1 & 2 & -1 & -1 & -1\\
1 & 2 & 1 & -1 & -1 & -1\\
2 & 1 & 1 & -1 & -1 & -1
\end{array}
\right)  ;
\end{array}
\]%
\[%
\begin{array}
[c]{ccc}%
M_{2}=\left(
\begin{array}
[c]{cccccc}%
1 & 1 & 1 & 1 & 1 & 2\\
1 & 1 & 1 & 1 & 2 & 1\\
1 & 1 & 1 & 2 & 2 & 2\\
1 & 1 & 2 & 2 & 1 & 2\\
1 & 2 & 2 & 1 & 1 & 2\\
2 & 1 & 2 & 2 & 2 & 2
\end{array}
\right)  , & M_{2}^{-1}= & \left(
\begin{array}
[c]{cccccc}%
-1 & -2 & 1 & -2 & 1 & 2\\
-2 & -2 & 2 & -2 & 2 & 1\\
1 & 2 & -2 & 2 & -1 & -1\\
-2 & -2 & 2 & -1 & 1 & 1\\
1 & 2 & -1 & 1 & -1 & -1\\
2 & 1 & -1 & 1 & -1 & -1
\end{array}
\right)  ;
\end{array}
\]%
\[%
\begin{array}
[c]{ccc}%
M_{3}=\left(
\begin{array}
[c]{cccccc}%
1 & 1 & 1 & 1 & 1 & 2\\
1 & 1 & 1 & 2 & 2 & 2\\
1 & 1 & 2 & 1 & 2 & 2\\
1 & 2 & 1 & 2 & 2 & 1\\
1 & 2 & 2 & 2 & 2 & 2\\
2 & 2 & 2 & 1 & 2 & 2
\end{array}
\right)  , & M_{3}^{-1}= & \left(
\begin{array}
[c]{cccccc}%
-2 & 2 & -2 & -2 & 1 & 2\\
2 & -2 & 1 & 2 & -1 & -1\\
-2 & 1 & -1 & -2 & 2 & 1\\
-2 & 2 & -2 & -2 & 2 & 1\\
1 & -1 & 2 & 2 & -2 & -1\\
2 & -1 & 1 & 1 & -1 & -1
\end{array}
\right)  ;
\end{array}
\]%
\[%
\begin{array}
[c]{ccc}%
M_{4}=\left(
\begin{array}
[c]{cccccc}%
1 & 1 & 1 & 2 & 2 & 2\\
1 & 2 & 2 & 1 & 1 & 2\\
1 & 2 & 2 & 2 & 2 & 2\\
2 & 1 & 2 & 1 & 2 & 1\\
2 & 1 & 2 & 2 & 2 & 2\\
2 & 2 & 2 & 1 & 2 & 2
\end{array}
\right)  , &  & M_{4}^{-1}=\left(
\begin{array}
[c]{cccccc}%
-2 & -2 & 1 & -2 & 2 & 2\\
-2 & -2 & 2 & -2 & 1 & 2\\
1 & 2 & -1 & 2 & -1 & -2\\
-2 & -2 & 2 & -2 & 2 & 1\\
2 & 1 & -1 & 2 & -2 & -1\\
2 & 2 & -2 & 1 & -1 & -1
\end{array}
\right)  ;
\end{array}
\]
or among $199$ special matrices containing at least one negative entry, e.g.,%

\[%
\begin{array}
[c]{ccc}%
N_{5}=\left(
\begin{array}
[c]{cccccc}%
1 & 1 & 1 & 1 & 1 & 1\\
1 & 1 & 1 & 1 & 2 & 2\\
1 & 1 & 1 & 2 & 1 & 2\\
1 & 1 & 2 & 1 & -1 & 1\\
1 & 2 & 1 & 2 & 2 & 2\\
1 & 2 & 2 & 1 & 2 & 2
\end{array}
\right)  , &  & N_{5}^{-1}=\left(
\begin{array}
[c]{cccccc}%
1 & 2 & -2 & 1 & 1 & -2\\
-1 & 1 & -2 & 1 & 2 & -1\\
1 & -2 & 2 & -1 & -2 & 2\\
1 & -2 & 2 & -1 & -1 & 1\\
1 & -1 & 1 & -1 & -1 & 1\\
-2 & 2 & -1 & 1 & 1 & -1
\end{array}
\right)  ;
\end{array}
\]%
\[%
\begin{array}
[c]{ccc}%
N_{6}=\left(
\begin{array}
[c]{cccccc}%
1 & 1 & 1 & 1 & 1 & 1\\
1 & 1 & 1 & 1 & 2 & 2\\
1 & 1 & 1 & 2 & 1 & 2\\
1 & 1 & 2 & 1 & 2 & 1\\
1 & 2 & 1 & 2 & 1 & 1\\
1 & 2 & -1 & 1 & 1 & 1
\end{array}
\right)  , &  & N_{6}^{-1}=\left(
\begin{array}
[c]{cccccc}%
2 & -2 & 1 & 1 & -2 & 1\\
1 & 2 & -2 & -2 & 2 & -1\\
1 & 1 & -1 & -1 & 1 & -1\\
-2 & -2 & 2 & 2 & -1 & 1\\
-2 & -1 & 1 & 2 & -1 & 1\\
1 & 2 & -1 & -2 & 1 & -1
\end{array}
\right)  ,
\end{array}
\]%
\[%
\begin{array}
[c]{ccccc}%
N_{7}=\left(
\begin{array}
[c]{cccccc}%
1 & 1 & 1 & 1 & 1 & 1\\
1 & 1 & 1 & 1 & 2 & 2\\
1 & 1 & 1 & 2 & 1 & 2\\
1 & 1 & 2 & 1 & 2 & 1\\
1 & 2 & 1 & 2 & 1 & 1\\
1 & 2 & 2 & 1 & 1 & -2
\end{array}
\right)  , &  & N_{7}^{-1}=\left(
\begin{array}
[c]{cccccc}%
2 & 1 & 1 & -2 & -2 & 1\\
1 & -1 & -2 & 1 & 2 & -1\\
1 & -2 & -1 & 2 & 1 & -1\\
-2 & 1 & 2 & -1 & -1 & 1\\
-2 & 2 & 1 & -1 & -1 & 1\\
1 & -1 & -1 & 1 & 1 & -1
\end{array}
\right)  , &  & \ldots,
\end{array}
\]%
\[%
\begin{array}
[c]{ccc}%
N_{203}=\left(
\begin{array}
[c]{cccccc}%
1 & 1 & 1 & 2 & 2 & 2\\
1 & 1 & 2 & 1 & 2 & 2\\
1 & -2 & 1 & -2 & -1 & -2\\
2 & -1 & 2 & -2 & -1 & -1\\
2 & -2 & 1 & -2 & -2 & -2\\
2 & -2 & 2 & -2 & -1 & -2
\end{array}
\right)  , &  & N_{203}^{-1}=\left(
\begin{array}
[c]{cccccc}%
2 & -2 & 1 & 2 & -1 & -1\\
1 & -2 & -1 & 2 & -2 & 1\\
-2 & 2 & -2 & -2 & 1 & 2\\
-1 & 1 & -2 & -2 & 1 & 2\\
2 & -2 & 2 & 2 & -2 & -1\\
-1 & 2 & 1 & -1 & 2 & -2
\end{array}
\right)  .
\end{array}
\]
Just as for $4\times4$ matrices, $\left\Vert \left(  M\;M^{-1}\right)
\right\Vert =2$ is valid. \ A full table of the special matrices would be
possible; we omit discussion of $N$ for reasons of space.

\section{Integer $7\times7$ Matrices}

Unlike the scenario for $6\times6$ matrices, there are only $2$ (not $4$)
equivalence classes for minimal unimodular zerofree $7\times7$ matrices $M$
with positive entries only:%
\[%
\begin{array}
[c]{ccc}%
M_{1}=\left(
\begin{array}
[c]{ccccccc}%
1 & 1 & 1 & 1 & 1 & 2 & 2\\
1 & 1 & 1 & 1 & 2 & 1 & 2\\
1 & 1 & 1 & 1 & 2 & 2 & 1\\
1 & 1 & 1 & 2 & 2 & 2 & 2\\
1 & 1 & 2 & 1 & 2 & 2 & 2\\
1 & 2 & 1 & 1 & 2 & 2 & 2\\
2 & 1 & 1 & 1 & 2 & 2 & 2
\end{array}
\right)  , &  & M_{1}^{-1}=\left(
\begin{array}
[c]{ccccccc}%
-2 & -2 & -2 & 1 & 1 & 1 & 2\\
-2 & -2 & -2 & 1 & 1 & 2 & 1\\
-2 & -2 & -2 & 1 & 2 & 1 & 1\\
-2 & -2 & -2 & 2 & 1 & 1 & 1\\
1 & 2 & 2 & -1 & -1 & -1 & -1\\
2 & 1 & 2 & -1 & -1 & -1 & -1\\
2 & 2 & 1 & -1 & -1 & -1 & -1
\end{array}
\right)  ;
\end{array}
\]%
\[%
\begin{array}
[c]{ccc}%
M_{2}=\left(
\begin{array}
[c]{ccccccc}%
1 & 1 & 1 & 1 & 1 & 1 & 2\\
1 & 1 & 1 & 1 & 1 & 2 & 1\\
1 & 1 & 1 & 1 & 2 & 1 & 1\\
1 & 1 & 1 & 2 & 1 & 1 & 1\\
1 & 2 & 2 & 2 & 2 & 2 & 2\\
2 & 1 & 2 & 2 & 2 & 2 & 2\\
2 & 2 & 1 & 2 & 2 & 2 & 2
\end{array}
\right)  , &  & M_{2}^{-1}=\left(
\begin{array}
[c]{ccccccc}%
-2 & -2 & -2 & -2 & 1 & 2 & 2\\
-2 & -2 & -2 & -2 & 2 & 1 & 2\\
-2 & -2 & -2 & -2 & 2 & 2 & 1\\
1 & 1 & 1 & 2 & -1 & -1 & -1\\
1 & 1 & 2 & 1 & -1 & -1 & -1\\
1 & 2 & 1 & 1 & -1 & -1 & -1\\
2 & 1 & 1 & 1 & -1 & -1 & -1
\end{array}
\right)  .
\end{array}
\]
Unlike the scenario for $5\times5$ matrices, however, $\left\Vert \left(
M\;M^{-1}\right)  \right\Vert =2$ is valid. \ Notice that $M_{2}=M_{1}%
^{\prime}$, the transpose of $M_{1}$. \ We have not attempted to find
equivalence classes for minimal $7\times7$ matrices $N$ containing at least
one negative entry.

Unimodular zerofree $8\times8$ matrices indeed exist. \ The $9\times9$ case
remains open. \ 

\section{Canonical Representatives}

Signed permutations act by shuffling rows and columns while optionally
flipping their signs; therefore many matrices that appear distinct are
actually identical up to this symmetry. Matrix classification implies
collapsing each entire orbit of such moves to a single, well-chosen
representative. Once this is done, it becomes possible to detect when two
matrices encode the same object and to avoid redundant variants. \ Among the
various well-known linear-algebraic \textquotedblleft canonical
forms\textquotedblright, the one that is closest in spirit to
signed-permutation canonicalization is the Smith normal form. \ Both solve a
similar shape of problem: a matrix modulo a group acting on the left and
right. \ Their computational philosophies are the same: normalize away
symmetries until only the essential features remain.

At the end of Section 1, a statement that $2\times2$ matrices $A$ and
$\tilde{A}$ are not equivalent (under left/right signed permutations) was left
unproven. \ Likewise, the inequivalence of $4\times4$ matrices $M_{1}$ and
$N_{2}$, fundamental to Section 3, was implied yet not demonstrated.

The following Magma program offers verification:%

\[%
\begin{array}
[c]{l}%
\text{\texttt{Z := IntegerRing();}}\\
\\
\text{\texttt{SignedPermutationGroup := function(n)\ \ }}\\
\text{\texttt{\ \ \ M := MatrixRing(Z, n); }}\\
\text{\texttt{\ \ \ Perms := [ M!PermutationMatrix(Z, p) : p in
SymmetricGroup(n) ]; }}\\
\text{\texttt{\ \ \ Signs := [ M!DiagonalMatrix([ s[i] : i in [1..n] ]) }}\\
\text{\texttt{\ \ \ \ \ \ : s in CartesianPower(\{-1,1\}, n) }}\mathtt{];}\\
\text{\texttt{\ \ \ return [ D * P : D in Signs, P in Perms ]; }}\\
\text{\texttt{end function; }}\\
\\
\text{\texttt{SignedPermutationDoubleAction := function(n) }}\\
\text{\texttt{\ \ \ S := SignedPermutationGroup(n); }}\\
\text{\texttt{\ \ \ return [
$<$%
P, Q%
$>$
: P in S, Q in S ];}}\\
\text{\texttt{end function;}}\\
\\
\text{\texttt{Action := function(g, M) }}\\
\text{\texttt{ \ \ return g[1] * M * g[2];}}\\
\text{\texttt{end function;}}\\
\\
\text{\texttt{CanonicalizeMatrix := function(A) }}\\
\text{\texttt{ \ \ n := Nrows(A);}}\\
\text{\texttt{ \ \ G := SignedPermutationDoubleAction(n); }}\\
\text{\texttt{ \ \ orbit := [ Action(g, A) : g in G ]; }}\\
\text{\texttt{ \ \ return Minimum(orbit); }}\\
\text{\texttt{end function; }}%
\end{array}
\bigskip
\]

\noindent In words, $\operatorname*{CanonicalizeMatrix}(A)$ returns the
smallest matrix in the orbit of $A$, comparing matrices lexicographically
using Magma's built-in row-major total ordering. \ Because%
\begin{align*}
\operatorname*{CanonicalizeMatrix}(A)  &  =\left(
\begin{array}
[c]{cc}%
-4 & -2\\
-1 & 1
\end{array}
\right) \\
&  \neq\left(
\begin{array}
[c]{cc}%
-4 & -1\\
-2 & 1
\end{array}
\right)  =\operatorname*{CanonicalizeMatrix}(\tilde{A}),
\end{align*}
the matrices $A$ and $\tilde{A}$ cannot occupy the same orbit, i.e., cannot be
equivalent. \ Similar reasoning applies to $M_{1}$ and $N_{2}$ because%
\begin{align*}
\operatorname*{CanonicalizeMatrix}(M_{1})  &  =\left(
\begin{array}
[c]{cccc}%
-2 & -2 & -2 & -1\\
-2 & -2 & -1 & -2\\
-2 & -1 & -1 & -1\\
-1 & -2 & -1 & -1
\end{array}
\right)  \neq\\
\left(
\begin{array}
[c]{cccc}%
-2 & -2 & -2 & -1\\
-2 & -1 & -1 & -1\\
-2 & 1 & 2 & -2\\
-1 & 1 & 2 & -1
\end{array}
\right)   &  =\operatorname*{CanonicalizeMatrix}(N_{2}).
\end{align*}
A conceptual proof for the $2\times2$ case is outlined in the Appendix.

The following Mathematica program also offers verification, although the
outcome differs slightly:%
\[%
\begin{array}
[c]{l}%
\text{\texttt{S = S[\#] = Join @@
Map[Permutations@*DiagonalMatrix]@Tuples[\{-1,1\}, \#] \&;}}\\
\\
\text{\texttt{G = G[\#] = Flatten[Table[\{P, Q\}, \{P, S[\#]\}, \{Q, S[\#]\}],
1] \&;}}\\
\\
\text{\texttt{Action[\{P\_, Q\_\}, M\_] := P.M.Q;}}\\
\\
\text{\texttt{intKey[n\_Integer] := \{If[n
$>$
0, 0, 1], Abs[n]\};}}\\
\\
\text{\texttt{matrixKey[X\_] := intKey /@ Flatten[X];}}\\
\\
\text{\texttt{CanonicalizeMatrix[M\_] := Module[\{n, orbit\}, n = Length[M];}%
}\\
\text{\texttt{\ \ \ orbit = Action[\#, M] \& /@ G[n];}}\\
\text{\texttt{ \ \ First@SortBy[orbit, matrixKey]];}}%
\end{array}
\]

\noindent In words, $\operatorname*{CanonicalizeMatrix}[A]$ returns the
smallest matrix in the orbit of $A$, comparing matrices by their row-major
flattenings, using a custom integer key that emulates the structural ordering%
\[
1<2<3<4<5<\cdots<-1<-2<-3<-4<-5<\cdots
\]
rather than Mathematica's native numerical ordering. \ Because
\begin{align*}
\operatorname*{CanonicalizeMatrix}[A]  &  =\left(
\begin{array}
[c]{cc}%
1 & 1\\
2 & -4
\end{array}
\right) \\
&  \neq\left(
\begin{array}
[c]{cc}%
1 & 2\\
1 & -4
\end{array}
\right)  =\operatorname*{CanonicalizeMatrix}[\tilde{A}],
\end{align*}
the matrices $A$ and $\tilde{A}$ cannot occupy the same orbit. \ Similar
reasoning applies to $M_{1}$ and $N_{2}$ because
\[
\operatorname*{CanonicalizeMatrix}[M_{1}]=M_{1}\neq N_{2}%
=\operatorname*{CanonicalizeMatrix}[N_{2}].
\]
This interesting result -- that the canonical representative of $M_{1}$
coincides with $M_{1}$ -- also applies to the $3\times3$ matrix $M$ in Section
2. \ 

Warp-cooperative CUDA becomes essential at the $6\times6$ scale because the
search space expands so precipitously that no single thread can traverse the
full signed-permutation orbit before the computation becomes impractically
slow. \ For $4\times4$, the group is small enough that a CPU can exhaustively
assess every candidate; but by $6\times6$, the orbit has grown to billions of
elements, and each candidate requires two full matrix multiplications and a
comparison. At that point, only a coordinated warp-level sweep -- where all
$32$ lanes advance through the orbit together and maintain a shared best
candidate -- keeps the computation efficient enough for the GPU to handle.

The sequence%
\[%
\begin{array}
[c]{ccccccccc}%
1, &  & \left(
\begin{array}
[c]{cc}%
1 & 1\\
1 & 2
\end{array}
\right)  , &  & \left(
\begin{array}
[c]{ccc}%
1 & 2 & 2\\
2 & 1 & 2\\
2 & 2 & 3
\end{array}
\right)  , &  & \left(
\begin{array}
[c]{cccc}%
1 & 1 & 1 & 2\\
1 & 1 & 2 & 1\\
1 & 2 & 2 & 2\\
2 & 1 & 2 & 2
\end{array}
\right)  , &  & \ldots
\end{array}
\]
captures our imagination. \ A question raised in \cite{Fnch-mtrx} -- what is
the next term in this list? -- has now been answered (in Section 4):%
\[
\left(
\begin{array}
[c]{ccccc}%
1 & 1 & 1 & 2 & 2\\
1 & 1 & 2 & 1 & 2\\
1 & 2 & 2 & 2 & 2\\
2 & 1 & 2 & 2 & 2\\
2 & 2 & 2 & 2 & 1
\end{array}
\right)
\]
but the term succeeding this is disappointingly non-unique (matrices
$M_{\iota}$ for $1\leq\iota\leq4$ in Section 5). \ No pattern is apparent.

\section{Acknowledgments}

The creators of Mathematica earn my gratitude every day:\ this paper could not
have otherwise been written. \ Jan Mangaldan's function
\textit{RandomUnimodularMatrix} (an implementation of a method given
in\ \cite{Hau1-mtrx, Hau2-mtrx}) was very helpful to me. \ For the first time,
I\ have used the Microsoft Copilot generative AI chatbot for assistance in
writing/testing code. \ It introduced me to the extensive \&\ very pertinent
software package Magma and to the extraordinary accelerated performance
associated with CUDA programming on Nvidia GPUs.

\section{Appendix}

Of the $32$ matrices in the orbit of $A$, exactly $4$ possess $-4$ in the
upper-left corner:%
\[%
\begin{array}
[c]{ccccccc}%
\left(
\begin{array}
[c]{cc}%
-4 & -2\\
-1 & 1
\end{array}
\right)  , &  & \left(
\begin{array}
[c]{cc}%
-4 & -2\\
1 & -1
\end{array}
\right)  , &  & \left(
\begin{array}
[c]{cc}%
-4 & 2\\
-1 & -1
\end{array}
\right)  , &  & \left(
\begin{array}
[c]{cc}%
-4 & 2\\
1 & 1
\end{array}
\right)  .
\end{array}
\]
Magna's row-major ordering compares first rows lexicographically, then second
rows. \ $\left(
\begin{array}
[c]{cc}%
-4 & -2
\end{array}
\right)  <\left(
\begin{array}
[c]{cc}%
-4 & 2
\end{array}
\right)  $; thus we keep%
\[%
\begin{array}
[c]{ccc}%
\left(
\begin{array}
[c]{cc}%
-4 & -2\\
-1 & 1
\end{array}
\right)  , &  & \left(
\begin{array}
[c]{cc}%
-4 & -2\\
1 & -1
\end{array}
\right)  .
\end{array}
\]
$\left(
\begin{array}
[c]{cc}%
-1 & 1
\end{array}
\right)  <\left(
\begin{array}
[c]{cc}%
1 & -1
\end{array}
\right)  $; thus we keep $\tbinom{-4\;-2}{-1\;\;1}\bigskip$.

Of the $32$ matrices in the orbit of $\tilde{A}$, exactly $4$ possess $-4$ in
the upper-left corner:%
\[%
\begin{array}
[c]{ccccccc}%
\left(
\begin{array}
[c]{cc}%
-4 & -1\\
-2 & 1
\end{array}
\right)  , &  & \left(
\begin{array}
[c]{cc}%
-4 & -1\\
2 & -1
\end{array}
\right)  , &  & \left(
\begin{array}
[c]{cc}%
-4 & 1\\
-2 & -1
\end{array}
\right)  , &  & \left(
\begin{array}
[c]{cc}%
-4 & 1\\
2 & 1
\end{array}
\right)  .
\end{array}
\]
\ $\left(
\begin{array}
[c]{cc}%
-4 & -1
\end{array}
\right)  <\left(
\begin{array}
[c]{cc}%
-4 & 1
\end{array}
\right)  $; thus we keep%
\[%
\begin{array}
[c]{ccc}%
\left(
\begin{array}
[c]{cc}%
-4 & -1\\
-2 & 1
\end{array}
\right)  , &  & \left(
\begin{array}
[c]{cc}%
-4 & -1\\
2 & -1
\end{array}
\right)  .
\end{array}
\]
$\left(
\begin{array}
[c]{cc}%
-2 & 1
\end{array}
\right)  <\left(
\begin{array}
[c]{cc}%
2 & -1
\end{array}
\right)  $; thus we keep $\tbinom{-4\;-1}{-2\;\;1}\bigskip$.

Of the $32$ matrices in the orbit of $A$, exactly $8$ possess $1$ in the
upper-left corner:%
\[%
\begin{array}
[c]{ccccccc}%
\left(
\begin{array}
[c]{cc}%
1 & 1\\
2 & -4
\end{array}
\right)  , &  & \left(
\begin{array}
[c]{cc}%
1 & 1\\
4 & -2
\end{array}
\right)  , &  & \left(
\begin{array}
[c]{cc}%
1 & 1\\
-2 & 4
\end{array}
\right)  , &  & \left(
\begin{array}
[c]{cc}%
1 & 1\\
-4 & 2
\end{array}
\right)  ,
\end{array}
\]%
\[%
\begin{array}
[c]{ccccccc}%
\left(
\begin{array}
[c]{cc}%
1 & -1\\
2 & 4
\end{array}
\right)  , &  & \left(
\begin{array}
[c]{cc}%
1 & -1\\
4 & 2
\end{array}
\right)  , &  & \left(
\begin{array}
[c]{cc}%
1 & -1\\
-2 & -4
\end{array}
\right)  , &  & \left(
\begin{array}
[c]{cc}%
1 & -1\\
-4 & -2
\end{array}
\right)  .
\end{array}
\]
Mathematica's ordering compares row-major flattenings of matrices, checking
entries one after another. \ The smallest possible second entry is $1$; thus
only
\[%
\begin{array}
[c]{ccccccc}%
\left(
\begin{array}
[c]{cc}%
1 & 1\\
2 & -4
\end{array}
\right)  , &  & \left(
\begin{array}
[c]{cc}%
1 & 1\\
4 & -2
\end{array}
\right)  , &  & \left(
\begin{array}
[c]{cc}%
1 & 1\\
-2 & 4
\end{array}
\right)  , &  & \left(
\begin{array}
[c]{cc}%
1 & 1\\
-4 & 2
\end{array}
\right)
\end{array}
\]
survive. \ The smallest possible third entry is $2$; thus only $\tbinom
{\,\,1\;\;\,1}{\,\,2\;-4}$ survives. \ There is no need to examine the fourth entry.

Of the $32$ matrices in the orbit of $\tilde{A}$, exactly $8$ possess $1$ in
the upper-left corner:%
\[%
\begin{array}
[c]{ccccccc}%
\left(
\begin{array}
[c]{cc}%
1 & 2\\
1 & -4
\end{array}
\right)  , &  & \left(
\begin{array}
[c]{cc}%
1 & 2\\
-1 & 4
\end{array}
\right)  , &  & \left(
\begin{array}
[c]{cc}%
1 & 4\\
1 & -2
\end{array}
\right)  , &  & \left(
\begin{array}
[c]{cc}%
1 & 4\\
-1 & 2
\end{array}
\right)  ,
\end{array}
\]%
\[%
\begin{array}
[c]{ccccccc}%
\left(
\begin{array}
[c]{cc}%
1 & -2\\
1 & 4
\end{array}
\right)  , &  & \left(
\begin{array}
[c]{cc}%
1 & -2\\
-1 & -4
\end{array}
\right)  , &  & \left(
\begin{array}
[c]{cc}%
1 & -4\\
1 & 2
\end{array}
\right)  , &  & \left(
\begin{array}
[c]{cc}%
1 & -4\\
-1 & -2
\end{array}
\right)  .
\end{array}
\]
The smallest possible second entry is $2$; thus only%
\[%
\begin{array}
[c]{ccc}%
\left(
\begin{array}
[c]{cc}%
1 & 2\\
1 & -4
\end{array}
\right)  , &  & \left(
\begin{array}
[c]{cc}%
1 & 2\\
-1 & 4
\end{array}
\right)
\end{array}
\]
survive. \ The smallest possible third entry is $1$; thus only $\tbinom
{\,\,1\;\;\,2}{\,\,1\;-4}$ survives. \ Again, there is no need to examine the
fourth entry.

\section{Addendum I}

Choosing a matrix $X$ uniformly at random from the pool of $3686400$ minimal
unimodular zerofree $5\times5$ matrices and defining $\Xi
=\operatorname*{CanonicalizeMatrix}[X]$, we obtain exact probabilities%
\[
\mathbb{P}(\Xi=M_{1})=\dfrac{1}{2}=\mathbb{P}(\Xi=N_{2}).
\]
The $6\times6$ analog requires more bookkeeping since there are $203$
equivalence classes to manage (see Addendum II). \pagebreak

Here are two $8\times8$ minimal unimodular zerofree matrices and their inverses:%

\[%
\begin{array}
[c]{ccc}%
\left(
\begin{array}
[c]{cccccccc}%
1 & 1 & 1 & 1 & 1 & 1 & 1 & 2\\
1 & 1 & 1 & 1 & 1 & 1 & 2 & 1\\
1 & 1 & 1 & 1 & 1 & 2 & 1 & 1\\
1 & 1 & 1 & 1 & 2 & 1 & 2 & 2\\
1 & 1 & 1 & 2 & 1 & 1 & 2 & 2\\
1 & 1 & 2 & 1 & 2 & 1 & 1 & 2\\
1 & 2 & 2 & 1 & 1 & 2 & 1 & 1\\
2 & 1 & 2 & 1 & 2 & 2 & 1 & 2
\end{array}
\right)  , &  & \left(
\begin{array}
[c]{cccccccc}%
-1 & -1 & -2 & 1 & 1 & -2 & 1 & 2\\
-1 & -2 & -2 & 2 & 1 & -2 & 2 & 1\\
1 & 2 & 1 & -2 & -1 & 2 & -1 & -1\\
-2 & -2 & -1 & 1 & 2 & -1 & 1 & 1\\
-2 & -2 & -1 & 2 & 1 & -1 & 1 & 1\\
1 & 1 & 2 & -1 & -1 & 1 & -1 & -1\\
1 & 2 & 1 & -1 & -1 & 1 & -1 & -1\\
2 & 1 & 1 & -1 & -1 & 1 & -1 & -1
\end{array}
\right)  ;
\end{array}
\]%
\[%
\begin{array}
[c]{ccc}%
\left(
\begin{array}
[c]{cccccccc}%
1 & 1 & 1 & 1 & 1 & 1 & 1 & 2\\
1 & 1 & 1 & 1 & 1 & 1 & 2 & 1\\
1 & 1 & 1 & 1 & 1 & 2 & 1 & 1\\
1 & 1 & 1 & 1 & 2 & 1 & 2 & 2\\
1 & 1 & 1 & 2 & 1 & 2 & 1 & 2\\
1 & 1 & 2 & 1 & 2 & 1 & 2 & 1\\
1 & 2 & 2 & 1 & 1 & 1 & 2 & 1\\
2 & 1 & 2 & 1 & 2 & 2 & 2 & 1
\end{array}
\right)  , &  & \left(
\begin{array}
[c]{cccccccc}%
-1 & -1 & -2 & 1 & 1 & -2 & 1 & 2\\
-2 & -2 & -1 & 2 & 1 & -2 & 2 & 1\\
2 & 1 & 1 & -2 & -1 & 2 & -1 & -1\\
-2 & -1 & -2 & 1 & 2 & -1 & 1 & 1\\
-2 & -2 & -1 & 2 & 1 & -1 & 1 & 1\\
1 & 1 & 2 & -1 & -1 & 1 & -1 & -1\\
1 & 2 & 1 & -1 & -1 & 1 & -1 & -1\\
2 & 1 & 1 & -1 & -1 & 1 & -1 & -1
\end{array}
\right)  .
\end{array}
\]
Just as for $4\times4$, $6\times6$ and $7\times7$ matrices, $\left\Vert
\left(  M\;M^{-1}\right)  \right\Vert =2$ is valid. \ Here also is one
$9\times9$ minimal unimodular zerofree matrix and its inverse:%

\[%
\begin{array}
[c]{ccc}%
\left(
\begin{array}
[c]{ccccccccc}%
1 & 1 & 1 & 1 & 1 & 1 & 1 & 1 & 2\\
1 & 1 & 1 & 1 & 1 & 1 & 1 & 2 & 1\\
1 & 1 & 1 & 1 & 1 & 2 & 2 & 1 & 2\\
1 & 1 & 1 & 1 & 2 & 1 & 2 & 2 & 1\\
1 & 1 & 2 & 2 & 1 & 1 & 1 & 1 & 2\\
1 & 2 & 1 & 2 & 2 & 2 & 1 & 2 & 2\\
1 & 2 & 2 & 2 & 2 & 2 & 2 & 2 & 2\\
2 & 2 & 1 & 1 & 1 & 1 & 1 & 2 & 1\\
2 & 2 & 1 & 2 & 2 & 2 & 2 & 2 & 2
\end{array}
\right)  , &  & \left(
\begin{array}
[c]{ccccccccc}%
-2 & -2 & 2 & 2 & 2 & 2 & -3 & 2 & -2\\
2 & 1 & -2 & -2 & -2 & -2 & 3 & -1 & 2\\
-2 & -2 & 2 & 2 & 2 & 2 & -2 & 2 & -3\\
1 & 2 & -2 & -2 & -1 & -2 & 2 & -2 & 3\\
-2 & -3 & 2 & 3 & 2 & 3 & -3 & 2 & -3\\
-3 & -2 & 3 & 2 & 2 & 3 & -3 & 2 & -3\\
2 & 2 & -2 & -2 & -2 & -3 & 3 & -2 & 3\\
1 & 2 & -1 & -1 & -1 & -1 & 1 & -1 & 1\\
2 & 1 & -1 & -1 & -1 & -1 & 1 & -1 & 1
\end{array}
\right)  .
\end{array}
\]
Just as for $5\times5$ matrices, $\left\Vert \left(  M\;M^{-1}\right)
\right\Vert =3$ is valid\ although $\left\Vert M\right\Vert =2<\left\Vert
M^{-1}\right\Vert $. \ The $10\times10$ case remains open.

\section{Addendum II}

The $203$ canonical representative matrices listed in this section are given
as vectors in row-major order.\pagebreak

\begin{center}%
\[%
\begin{tabular}
[c]{|l|}\hline
Matrices $M_{1},M_{2},M_{3},M_{4}$ containing no negative entries:\\\hline
{\small 1, 1, 1, 1, 1, 2, 1, 1, 1, 1, 2, 1, 1, 1, 1, 2, 1, 1, 1, 1, 2, 1, 2,
2, 1, 2, 1, 2, 1, 2, 2, 1, 1, 2, 2, 1}\\\hline
{\small 1, 1, 1, 1, 1, 2, 1, 1, 1, 1, 2, 1, 1, 1, 1, 2, 2, 2, 1, 1, 2, 2, 1,
2, 1, 2, 2, 1, 1, 2, 2, 1, 2, 2, 2, 2}\\\hline
{\small 1, 1, 1, 1, 1, 2, 1, 1, 1, 2, 2, 2, 1, 1, 2, 1, 2, 2, 1, 2, 1, 2, 2,
1, 1, 2, 2, 2, 2, 2, 2, 2, 2, 1, 2, 2}\\\hline
{\small 1, 1, 1, 2, 2, 2, 1, 2, 2, 1, 1, 2, 1, 2, 2, 2, 2, 2, 2, 1, 2, 1, 2,
1, 2, 1, 2, 2, 2, 2, 2, 2, 2, 1, 2, 2}\\\hline
\end{tabular}
\ \ \ \
\]

\end{center}

\noindent Choosing a matrix $X$ uniformly at random from the pool of
$196942233600$ minimal unimodular zerofree $6\times6$ matrices and defining
$\Xi=\operatorname*{CanonicalizeMatrix}[X]$, we obtain exact probabilities%
\[%
\begin{array}
[c]{ccc}%
\mathbb{P}(\Xi=M_{1})=\mathbb{P}(\Xi=M_{4})=\dfrac{1}{1113}, &  &
\mathbb{P}(\Xi=M_{2})=\mathbb{P}(\Xi=M_{3})=\dfrac{2}{371}.
\end{array}
\]
The odds of $M_{1}$ (or of $M_{2}$) occurring are the smallest among all $203$
ratios appearing here; $1/1113$ appears with $N_{196}$ as well.%
\[%
\begin{tabular}
[c]{|l|}\hline
Matrices $N_{5},N_{6},\ldots,N_{17}$ containing exactly $1$ negative
entry:\\\hline
{\small 1, 1, 1, 1, 1, 1, 1, 1, 1, 1, 2, 2, 1, 1, 1, 2, 1, 2, 1, 1, 2, 1, -1,
1, 1, 2, 1, 2, 2, 2, 1, 2, 2, 1, 2, 2}\\\hline
{\small 1, 1, 1, 1, 1, 1, 1, 1, 1, 1, 2, 2, 1, 1, 1, 2, 1, 2, 1, 1, 2, 1, 2,
1, 1, 2, 1, 2, 1, 1, 1, 2, -1, 1, 1, 1}\\\hline
{\small 1, 1, 1, 1, 1, 1, 1, 1, 1, 1, 2, 2, 1, 1, 1, 2, 1, 2, 1, 1, 2, 1, 2,
1, 1, 2, 1, 2, 1, 1, 1, 2, 2, 1, 1, -2}\\\hline
{\small 1, 1, 1, 1, 1, 1, 1, 1, 1, 1, 2, 2, 1, 1, 1, 2, 1, 2, 1, 1, 2, 1, 2,
1, 1, 2, 1, 2, 2, 2, 1, 2, -1, 1, 2, 2}\\\hline
{\small 1, 1, 1, 1, 1, 1, 1, 1, 1, 1, 2, 2, 1, 1, 1, 2, 1, 2, 1, 1, 2, 1, 2,
1, 1, 2, 1, 2, 2, 2, 1, 2, 2, 1, 2, -1}\\\hline
{\small 1, 1, 1, 1, 1, 1, 1, 1, 1, 1, 2, 2, 1, 1, 1, 2, 1, 2, 1, 1, 2, 2, 2,
2, 1, 2, 1, 2, -2, 1, 1, 2, 2, 2, 1, 2}\\\hline
{\small 1, 1, 1, 1, 1, 1, 1, 1, 1, 1, 2, 2, 1, 1, 1, 2, 1, 2, 1, 1, 2, 2, 2,
2, 1, 2, 2, 1, 2, 2, 1, -1, 1, 1, 2, 1}\\\hline
{\small 1, 1, 1, 1, 1, 1, 1, 1, 1, 1, 2, 2, 1, 1, 1, 2, 1, 2, 1, 2, 2, 1, 1,
1, 2, 1, 2, 1, 2, 2, 2, -1, 1, 2, 1, 2}\\\hline
{\small 1, 1, 1, 1, 1, 1, 1, 1, 1, 1, 2, 2, 1, 1, 2, 2, 1, 1, 1, 2, 1, 2, 2,
2, 1, 2, 2, 2, 1, 2, 2, -1, 2, 1, 2, 1}\\\hline
{\small 1, 1, 1, 1, 1, 1, 1, 1, 1, 1, 2, 2, 1, 1, 2, 2, 2, 2, 1, 2, 1, 2, 2,
2, 2, 1, 2, 2, 1, 2, 2, 2, 1, 2, -1, 1}\\\hline
{\small 1, 1, 1, 1, 1, 2, 1, 1, 1, 2, 2, 2, 1, 1, 2, 1, 2, -1, 1, 1, 2, 2, 2,
1, 1, 2, 2, 1, 1, 2, 1, 2, 2, 1, 2, 1}\\\hline
{\small 1, 1, 1, 1, 2, 2, 1, 1, 1, 2, 1, 2, 1, 1, 2, 1, 2, 1, 1, 1, 2, 2, 2,
2, 1, 2, 2, 2, 1, 2, 1, 2, 2, 2, -1, 1}\\\hline
{\small 1, 1, 1, 1, 2, 2, 1, 1, 1, 2, 1, 2, 1, 2, 2, 1, 1, 1, 1, 2, 2, 1, 2,
2, 2, 1, 2, 1, 2, 2, 2, 1, 2, -1, 2, 1}\\\hline
\end{tabular}
\ \ \ \ \ \
\]

\noindent Here%
\begin{align*}
\mathbb{P}(\Xi &  =N_{5})=\mathbb{P}(\Xi=N_{6})=\mathbb{P}(\Xi=N_{8}%
)=\mathbb{P}(\Xi=N_{9})=\mathbb{P}(\Xi=N_{10})=\mathbb{P}(\Xi=N_{11})\\
&  =\mathbb{P}(\Xi=N_{12})=\mathbb{P}(\Xi=N_{14})=\mathbb{P}(\Xi
=N_{16})=\mathbb{P}(\Xi=N_{17})=\dfrac{2}{371},
\end{align*}%
\[
\mathbb{P}(\Xi=N_{7})=\mathbb{P}(\Xi=N_{13})=\mathbb{P}(\Xi=N_{15})=\dfrac
{1}{371}.
\]

\[%
\begin{tabular}
[c]{|l|}\hline
Matrices $N_{18},N_{19},\ldots,N_{41}$ containing exactly $3$ negative
entries:\\\hline
{\small 1, 1, 1, 1, 1, 1, 1, 1, 1, 1, 2, -1, 1, 1, 2, 2, 2, -1, 1, 2, 1, 2, 2,
-1, 2, 1, 2, 2, 1, 2, 2, 2, 1, 2, 2, 1}\\\hline
{\small 1, 1, 1, 1, 1, 1, 1, 1, 1, 1, 2, 2, 1, 1, 1, -1, -1, -2, 1, 1, 2, 2,
2, 2, 1, 2, 2, 1, 2, 2, 2, 1, 2, 2, 2, 1}\\\hline
{\small 1, 1, 1, 1, 1, 1, 1, 1, 1, 1, 2, 2, 1, 1, 1, -1, 1, 2, 1, 1, 2, 2, 2,
2, 1, 2, 1, -1, 1, 1, 1, 2, 2, -1, 1, 2}\\\hline
{\small 1, 1, 1, 1, 1, 1, 1, 1, 1, 1, 2, 2, 1, 1, 1, 2, 1, -1, 1, 1, 2, 2, 1,
-2, 1, 2, 1, 2, 2, 2, 1, 2, 2, 2, 1, -1}\\\hline
{\small 1, 1, 1, 1, 1, 1, 1, 1, 1, 1, 2, 2, 1, 1, 1, 2, 1, -1, 1, 1, 2, 2, 1,
1, 1, 2, 1, 2, 2, -1, 1, 2, 2, 2, 1, -1}\\\hline
{\small 1, 1, 1, 1, 1, 1, 1, 1, 1, 1, 2, 2, 1, 1, 1, 2, 1, 2, 1, 1, 2, -1, -2,
-2, 1, 2, 2, 1, 1, 1, 2, 1, 2, 2, 1, 2}\\\hline
{\small 1, 1, 1, 1, 1, 1, 1, 1, 1, 1, 2, 2, 1, 1, 1, 2, 1, 2, 1, 1, 2, -1, -1,
-1, 1, 2, 2, 1, 2, 2, 2, 1, 2, 2, 1, 2}\\\hline
{\small 1, 1, 1, 1, 1, 1, 1, 1, 1, 1, 2, 2, 1, 1, 1, 2, 1, 2, 1, 1, 2, 1, 2,
1, 1, 2, 1, 2, 2, 2, 1, -1, -1, 1, -1, 2}\\\hline
{\small 1, 1, 1, 1, 1, 1, 1, 1, 1, 1, 2, 2, 1, 1, 1, 2, 1, 2, 1, 1, 2, 1, 2,
1, 1, 2, 1, 2, 2, 2, 1, -1, 2, -2, 2, -1}\\\hline
{\small 1, 1, 1, 1, 1, 1, 1, 1, 1, 1, 2, 2, 1, 1, 1, 2, 1, 2, 1, 1, 2, 1, 2,
1, 1, 2, 1, 2, 2, 2, 1, -1, 2, 1, -1, -1}\\\hline
{\small 1, 1, 1, 1, 1, 1, 1, 1, 1, 1, 2, 2, 1, 1, 1, 2, 1, 2, 1, 1, 2, 2, 2,
2, 1, 2, 2, 1, 2, 2, 1, -1, -2, 1, -1, 1}\\\hline
{\small 1, 1, 1, 1, 1, 1, 1, 1, 1, 1, 2, 2, 1, 1, 1, 2, 1, 2, 1, 2, 2, 1, 1,
1, 1, -1, -2, 1, -1, 1, 2, 2, 1, 1, 2, 2}\\\hline
{\small 1, 1, 1, 1, 1, 1, 1, 1, 1, 1, 2, 2, 1, 1, 1, 2, 1, 2, 1, 2, 2, 1, 1,
1, 1, -1, -2, 2, -2, 1, 2, 2, 1, 2, 1, 2}\\\hline
{\small 1, 1, 1, 1, 1, 1, 1, 1, 1, 1, 2, 2, 1, 1, 2, 2, 2, 2, 1, -1, 1, 2, 1,
1, 1, -1, 2, 2, 1, 2, 1, -2, 1, 2, 1, 2}\\\hline
{\small 1, 1, 1, 1, 1, 1, 1, 1, 1, 1, 2, 2, 1, 1, 2, 2, 2, 2, 1, 2, 1, 2, 2,
2, 1, 2, -2, 1, -1, -2, 2, 2, 1, 2, 2, 1}\\\hline
{\small 1, 1, 1, 1, 1, 1, 1, 1, 1, 1, 2, 2, 1, 1, 2, 2, 2, 2, 1, 2, 1, 2, 2,
2, 2, 1, 2, 2, 1, 2, 2, -1, 1, -1, -1, 1}\\\hline
{\small 1, 1, 1, 1, 1, 2, 1, 1, 1, 1, 2, 1, 1, 1, 1, 2, 2, 2, 1, 1, -1, 2, 1,
2, 1, 2, -1, 1, 1, 2, 1, 2, -2, 1, 1, 1}\\\hline
{\small 1, 1, 1, 1, 1, 2, 1, 1, 1, 1, 2, 1, 1, 1, 1, 2, 2, 2, 1, 1, -1, 2, 1,
2, 1, 2, -1, 2, 2, 2, 1, 2, -2, 2, 2, 1}\\\hline
{\small 1, 1, 1, 1, 1, 2, 1, 1, 1, 1, 2, 1, 1, 1, 1, 2, 2, 2, 1, 1, 2, -1, 1,
2, 1, 2, 1, -2, 1, 1, 1, 2, 2, -2, 1, 2}\\\hline
{\small 1, 1, 1, 1, 1, 2, 1, 1, 1, 1, 2, 1, 1, 1, 1, 2, 2, 2, 1, 1, 2, -1, 1,
2, 1, 2, 1, -1, 2, 1, 1, 2, 2, -1, 2, 2}\\\hline
{\small 1, 1, 1, 1, 1, 2, 1, 1, 1, 2, 2, 2, 1, 1, 2, 1, 2, 2, 1, -1, 1, 2, 2,
1, 1, -1, 2, 2, 2, 2, 1, -2, 1, 2, 1, 1}\\\hline
{\small 1, 1, 1, 1, 1, 2, 1, 1, 1, 2, 2, 2, 1, 1, 2, 1, 2, 2, 1, -1, 1, 2, 2,
1, 2, -1, 2, 1, 2, 2, 2, -2, 1, 1, 1, 1}\\\hline
{\small 1, 1, 1, 1, 1, 2, 1, 1, 1, 2, 2, 2, 1, 1, 2, 1, 2, 2, 1, 1, 2, -2, 1,
1, 1, 2, 1, -1, 1, 2, 1, 2, 1, -2, 1, 1}\\\hline
{\small 1, 1, 1, 1, 1, 2, 1, 1, 1, 2, 2, 2, 1, 1, 2, 1, 2, 2, 1, 1, 2, -2, 1,
1, 1, 2, 2, -1, 2, 2, 1, 2, 2, -2, 2, 1}\\\hline
\end{tabular}
\ \
\]

\noindent Only $N_{24}$, $N_{37}$ have probabilities $1/371$; the other $N$s
have probabilities $2/371.$%

\[%
\begin{tabular}
[c]{|l|}\hline
Matrices $N_{42},N_{43},\ldots,N_{70}$ containing exactly $4$ negative
entries:\\\hline
{\small 1, 1, 1, 1, 1, 2, 1, 1, 1, 1, 2, 1, 1, 1, 1, 2, 2, 2, 1, 1, 2, -2, 1,
-2, 1, 2, 2, 1, 1, 2, 2, 1, 2, -2, 2, -2}\\\hline
{\small 1, 1, 1, 1, 1, 2, 1, 1, 1, 1, 2, 1, 1, 1, 1, 2, 2, 2, 1, 1, 2, -2, 1,
-2, 1, 2, 2, 1, 2, 1, 2, 1, 2, -2, 1, -1}\\\hline
{\small 1, 1, 1, 1, 1, 2, 1, 1, 1, 1, 2, 1, 1, 1, 1, 2, 2, 2, 1, 1, 2, -1, 1,
-1, 1, 2, 2, 1, 1, 2, 2, 1, 2, -1, 2, -1}\\\hline
{\small 1, 1, 1, 1, 1, 2, 1, 1, 1, 1, 2, 1, 1, 1, 1, 2, 2, 2, 1, 1, 2, -1, 1,
-1, 1, 2, 2, 2, 2, 2, 2, 1, 2, -2, 1, -1}\\\hline
{\small 1, 1, 1, 1, 1, 2, 1, 1, 1, 1, 2, 1, 1, 1, 1, 2, 2, 2, 1, 1, 2, -1, 2,
-2, 1, 2, 2, 1, 2, 1, 2, 1, 2, -1, 2, -1}\\\hline
{\small 1, 1, 1, 1, 1, 2, 1, 1, 1, 1, 2, 1, 1, 1, 1, 2, 2, 2, 1, 1, 2, -1, 2,
-2, 1, 2, 2, 2, 2, 2, 2, 1, 2, -2, 2, -2}\\\hline
{\small 1, 1, 1, 1, 1, 2, 1, 1, 1, 1, 2, 1, 1, 1, 1, 2, 2, 2, 1, 2, 2, 1, 1,
2, 1, -1, -2, 1, 1, 1, 2, -1, -2, 1, 2, 1}\\\hline
{\small 1, 1, 1, 1, 1, 2, 1, 1, 1, 1, 2, 1, 1, 1, 1, 2, 2, 2, 1, 2, 2, 1, 1,
2, 1, -1, -2, 2, 1, 2, 2, -1, -2, 2, 2, 2}\\\hline
{\small 1, 1, 1, 1, 1, 2, 1, 1, 1, 2, 2, 2, 1, 1, -1, 1, 2, -1, 1, 2, 2, 2, 2,
2, 1, 2, -1, 1, 2, -2, 2, 2, 1, 2, 2, 2}\\\hline
{\small 1, 1, 1, 1, 1, 2, 1, 1, 1, 2, 2, 2, 1, 1, 1, -1, -2, 1, 1, 1, 2, 1, 2,
2, 1, 2, 2, 1, 1, 2, 2, 1, 2, -1, -1, 2}\\\hline
{\small 1, 1, 1, 1, 1, 2, 1, 1, 1, 2, 2, 2, 1, 1, 1, -1, -2, 1, 1, 1, 2, 1, 2,
2, 1, 2, 2, 2, 2, 2, 2, 1, 2, -2, -2, 2}\\\hline
{\small 1, 1, 1, 1, 1, 2, 1, 1, 1, 2, 2, 2, 1, 1, 1, -1, -2, 1, 1, 1, 2, 2, 2,
1, 1, 2, 2, 2, 2, 2, 2, 1, 2, -1, -2, 1}\\\hline
{\small 1, 1, 1, 1, 1, 2, 1, 1, 1, 2, 2, 2, 1, 1, 2, 1, 2, 2, 1, 2, 1, -2, -2,
1, 1, 2, 2, -2, -2, 2, 2, 2, 2, 1, 2, 2}\\\hline
{\small 1, 1, 1, 1, 1, 2, 1, 1, 1, 2, 2, 2, 1, 1, 2, 1, 2, 2, 1, 2, 1, -2, 1,
-2, 1, 2, 2, -2, 1, -1, 2, 2, 2, 1, 2, 2}\\\hline
{\small 1, 1, 1, 1, 1, 2, 1, 1, 1, 2, 2, 2, 1, 1, 2, 1, 2, 2, 1, 2, 1, -1, -1,
1, 1, 2, 2, -1, -1, 2, 2, 2, 2, 1, 2, 2}\\\hline
{\small 1, 1, 1, 1, 1, 2, 1, 1, 1, 2, 2, 2, 1, 1, 2, 1, 2, 2, 1, 2, 1, -1, -1,
1, 1, 2, 2, 2, 2, 2, 2, 2, 2, -2, -1, 2}\\\hline
{\small 1, 1, 1, 1, 1, 2, 1, 1, 1, 2, 2, 2, 1, 1, 2, 1, 2, 2, 1, 2, 1, -1, 2,
-2, 1, 2, 2, -1, 2, -1, 2, 2, 2, 1, 2, 2}\\\hline
{\small 1, 1, 1, 1, 1, 2, 1, 1, 1, 2, 2, 2, 1, 1, 2, 1, 2, 2, 1, 2, 1, 2, 1,
2, 1, -1, -2, 1, 1, 1, 2, -1, -2, 2, 2, 2}\\\hline
{\small 1, 1, 1, 1, 1, 2, 1, 1, 1, 2, 2, 2, 1, 1, 2, 1, 2, 2, 1, 2, 1, 2, 1,
2, 1, 2, -2, 1, -2, 1, 2, 2, -2, 2, -1, 2}\\\hline
{\small 1, 1, 1, 1, 1, 2, 1, 1, 1, 2, 2, 2, 1, 1, 2, 1, 2, 2, 1, 2, 2, -1, 2,
-1, 1, 2, 2, -2, 2, -2, 2, 2, 1, 2, 2, 2}\\\hline
{\small 1, 1, 1, 1, 1, 2, 1, 1, 1, 2, 2, 2, 1, 1, 2, 1, 2, 2, 1, 2, 2, 2, 2,
2, 1, -2, 1, -1, 1, 1, 2, -2, 2, -2, 1, 2}\\\hline
{\small 1, 1, 1, 1, 1, 2, 1, 1, 1, 2, 2, 2, 1, 1, 2, 1, 2, 2, 1, 2, 2, 2, 2,
2, 1, -1, 1, -1, 2, 1, 2, -1, 2, -2, 2, 2}\\\hline
{\small 1, 1, 1, 1, 1, 2, 1, 1, 1, 2, 2, 2, 1, 1, 2, 2, 2, 1, 1, 2, 1, -1, -2,
2, 1, 2, 2, 2, 2, 2, 2, 2, 2, -1, -2, 2}\\\hline
{\small 1, 1, 1, 1, 1, 2, 1, 1, 1, 2, 2, 2, 1, 1, 2, 2, 2, 1, 1, 2, 2, 2, 2,
2, 1, -2, 1, 1, 2, -2, 2, -2, 2, 1, 2, -2}\\\hline
{\small 1, 1, 1, 1, 1, 2, 1, 1, 1, 2, 2, 2, 1, 1, 2, 2, 2, 1, 1, 2, 2, 2, 2,
2, 1, -1, 1, 1, 2, -1, 2, -1, 2, 1, 2, -1}\\\hline
{\small 1, 1, 1, 1, 1, 2, 1, 1, 1, 2, 2, 2, 1, 2, 2, 2, 2, 2, 2, 1, 2, 2, 2,
2, 2, -2, 1, 1, 2, -1, 2, -2, 2, 1, 2, -2}\\\hline
{\small 1, 1, 1, 2, 2, 2, 1, 1, 2, 1, 2, 2, 1, 1, 2, -2, 1, -2, 1, 2, 2, 1, 2,
1, 1, 2, 2, 2, 2, 2, 2, 2, 2, -1, 2, -2}\\\hline
{\small 1, 1, 1, 2, 2, 2, 1, 2, 2, 1, 1, 2, 1, 2, 2, 2, 2, 2, 2, 1, 2, 2, 2,
2, 2, -1, 2, 1, 2, -1, 2, -2, 2, 1, 2, -2}\\\hline
{\small 1, 1, 1, 2, 2, 2, 1, 2, 2, 1, 1, 2, 1, 2, 2, 2, 2, 2, 2, 1, 2, 2, 2,
2, 2, 1, 2, -1, -2, 1, 2, 2, 2, -1, -2, 2}\\\hline
\end{tabular}
\ \
\]

\noindent All $N$s have probabilities $2/371.$%

\[%
\begin{tabular}
[c]{|l|}\hline
Matrices $N_{71},N_{72},\ldots,N_{94}$ containing exactly $5$ negative
entries:\\\hline
{\small 1, 1, 1, 1, 1, 1, 1, 1, 1, 1, 2, 2, 1, 1, 1, -1, 1, 2, 1, 1, 2, 2, 2,
2, 1, 2, -1, -2, 2, 2, 1, 2, -2, -2, 2, 1}\\\hline
{\small 1, 1, 1, 1, 1, 1, 1, 1, 1, 1, 2, 2, 1, 1, 1, 2, -1, -2, 1, 1, 2, 2, 1,
1, 1, 2, 1, 2, 2, 2, 1, 2, -1, 2, -1, -2}\\\hline
{\small 1, 1, 1, 1, 1, 1, 1, 1, 1, 1, 2, 2, 1, 1, 1, 2, -1, -2, 1, 1, 2, 2, 2,
2, 1, 2, 2, 2, 2, 1, 1, -1, 1, 2, -2, -2}\\\hline
{\small 1, 1, 1, 1, 1, 1, 1, 1, 1, 1, 2, 2, 1, 1, 1, 2, 1, -1, 1, 1, 2, 2, 2,
2, 1, 2, -1, 1, 2, -1, 1, 2, -2, 1, 2, -2}\\\hline
{\small 1, 1, 1, 1, 1, 1, 1, 1, 1, 1, 2, 2, 1, 1, 1, 2, 1, 2, 1, 2, 2, 1, 2,
2, 1, -1, -2, -1, -1, -1, 2, 2, 1, 2, 1, 2}\\\hline
{\small 1, 1, 1, 1, 1, 1, 1, 1, 1, 1, 2, 2, 1, 1, 2, 2, 1, 1, 1, 2, 1, -1, 2,
2, 1, -1, -1, 2, 1, 2, 1, -2, -2, 2, 1, 2}\\\hline
{\small 1, 1, 1, 1, 1, 1, 1, 1, 1, 1, 2, 2, 1, 1, 2, 2, 2, 2, 1, -1, 1, 2, 1,
1, 2, -1, -1, 1, 1, 2, 2, -2, -2, 1, 1, 2}\\\hline
{\small 1, 1, 1, 1, 1, 1, 1, 1, 1, 1, 2, 2, 1, 1, 2, 2, 2, 2, 1, 1, -1, -2, 1,
2, 1, 2, -1, -2, 1, 1, 1, -1, 2, 2, 1, 2}\\\hline
{\small 1, 1, 1, 1, 1, 1, 1, 1, 1, 1, 2, 2, 1, 1, 2, 2, 2, 2, 1, 1, -1, -2, 1,
2, 1, 2, -1, -1, 1, 2, 1, -1, 2, 1, 1, 1}\\\hline
{\small 1, 1, 1, 1, 1, 1, 1, 1, 1, 1, 2, 2, 1, 1, 2, 2, 2, 2, 1, 1, -1, -2, 1,
2, 1, 2, 2, 2, 1, -1, 1, 2, -1, -2, 1, 1}\\\hline
{\small 1, 1, 1, 1, 1, 1, 1, 1, 1, 1, 2, 2, 1, 1, 2, 2, 2, 2, 1, 2, -1, -1, 1,
2, 1, 2, -1, -2, 1, 1, 2, 2, 1, -1, 2, 1}\\\hline
{\small 1, 1, 1, 1, 1, 1, 1, 1, 1, 1, 2, 2, 1, 1, 2, 2, 2, 2, 1, 2, 1, -1, 1,
1, 1, -1, -1, 2, 1, 2, 1, -2, -2, 2, 1, 2}\\\hline
{\small 1, 1, 1, 1, 1, 1, 1, 1, 1, 1, 2, 2, 1, 1, 2, 2, 2, 2, 1, 2, 1, -1, 1,
1, 2, 1, -2, -2, 1, 2, 2, 2, -1, -2, 1, 2}\\\hline
{\small 1, 1, 1, 1, 1, 1, 1, 1, 1, 1, 2, 2, 1, 1, 2, 2, 2, 2, 1, 2, 1, 2, 2,
2, 1, -1, -2, -2, -1, -2, 2, 2, 1, 2, 2, 1}\\\hline
{\small 1, 1, 1, 1, 1, 1, 1, 1, 1, 1, 2, 2, 1, 1, 2, 2, 2, 2, 1, 2, 2, -1, 1,
-1, 1, 2, 2, -2, 1, -2, 2, 2, 1, -1, 2, 1}\\\hline
{\small 1, 1, 1, 1, 1, 2, 1, 1, 1, 2, -2, 1, 1, 2, 2, 1, -1, 1, 1, 2, 2, 2,
-2, 1, 2, 1, 2, 2, -1, 2, 2, 1, 2, 2, -2, 1}\\\hline
{\small 1, 1, 1, 1, 1, 2, 1, 1, 1, 2, -1, 2, 1, 1, 1, 2, -2, 1, 1, 1, 2, 1,
-2, 1, 1, 2, 1, 2, -1, 1, 1, 2, 2, 2, -2, 1}\\\hline
{\small 1, 1, 1, 1, 1, 2, 1, 1, 1, 2, -1, 2, 1, 1, 1, 2, -2, 1, 1, 1, 2, 1,
-1, 2, 1, 2, 1, 2, -1, 1, 1, 2, 2, 2, -1, 2}\\\hline
{\small 1, 1, 1, 1, 1, 2, 1, 1, 1, 2, -1, 2, 1, 1, 2, 1, -2, 1, 1, 1, 2, 2,
-2, 2, 1, 2, 2, 2, -1, 2, 1, 2, 2, 2, -2, 1}\\\hline
{\small 1, 1, 1, 1, 1, 2, 1, 1, 1, 2, -1, 2, 1, 2, 2, 1, -1, 1, 1, 2, 2, 2,
-1, 2, 2, 1, 2, 2, -1, 2, 2, 1, 2, 2, -2, 1}\\\hline
{\small 1, 1, 1, 1, 1, 2, 1, 1, 1, 2, 2, 2, 1, 1, 2, 2, 2, 1, 1, 1, 2, -1, -2,
-1, 1, 2, 2, 1, 1, 2, 1, 2, 2, -1, -2, 1}\\\hline
{\small 1, 1, 1, 1, 1, 2, 1, 1, 1, 2, 2, 2, 1, 1, 2, 2, 2, 1, 1, 2, 2, 1, 1,
2, 1, -1, -1, 1, 2, 1, 1, -2, -1, 1, 2, -1}\\\hline
{\small 1, 1, 1, 1, 2, 2, 1, 1, 1, 2, 1, 2, 1, 1, 2, 1, -2, -1, 1, 1, 2, -1,
-2, -2, 1, 2, 2, 1, 1, 1, 1, 2, 2, 1, 2, 2}\\\hline
{\small 1, 1, 1, 1, 2, 2, 1, 2, 2, -1, 1, 2, 1, 2, 2, -2, 1, 1, 2, 1, 2, -1,
2, 1, 2, 1, 2, -2, 1, 1, 2, 2, 2, -1, 2, 2}\\\hline
\end{tabular}
\
\]

\noindent Only $N_{94}$ has probability $1/371$; the other $N$s have
probabilities $2/371.$%

\[%
\begin{tabular}
[c]{|l|}\hline
Matrices $N_{95},N_{96},\ldots,N_{111}$ containing exactly $7$ negative
entries:\\\hline
{\small 1, 1, 1, 1, 1, 1, 1, 1, 1, -1, -1, -2, 1, 1, 2, 2, 2, 2, 1, 2, -1, 2,
-1, 1, 1, 2, -2, 2, -2, 1, 2, 2, 2, 2, 1, 1}\\\hline
{\small 1, 1, 1, 1, 1, 1, 1, 1, 1, 1, 2, 2, 1, 1, 1, -1, -1, -2, 1, 2, 2, 1,
2, 2, 1, -1, -2, 2, 2, 2, 2, -1, -2, 2, 2, 1}\\\hline
{\small 1, 1, 1, 1, 1, 1, 1, 1, 1, 1, 2, 2, 1, 1, 2, 2, 1, 1, 1, 2, 1, 2, 2,
2, 1, -2, 1, -1, -1, -2, 2, -2, 1, -2, 1, -1}\\\hline
{\small 1, 1, 1, 1, 1, 1, 1, 1, 1, 1, 2, 2, 1, 1, 2, 2, 2, 2, 1, -1, 1, 2, -2,
-2, 1, -1, 2, 2, 1, 2, 1, -2, 1, 2, -2, -1}\\\hline
{\small 1, 1, 1, 1, 1, 1, 1, 1, 1, 1, 2, 2, 1, 1, 2, 2, 2, 2, 1, -1, 1, 2, 1,
1, 1, -1, 2, 2, -1, -2, 1, -2, 1, 2, -1, -2}\\\hline
{\small 1, 1, 1, 1, 1, 1, 1, 1, 1, 1, 2, 2, 1, 1, 2, 2, 2, 2, 1, 1, -1, -2, 1,
-1, 1, 2, -2, -1, 2, 2, 1, 2, -2, -2, 2, 1}\\\hline
{\small 1, 1, 1, 1, 1, 1, 1, 1, 1, 1, 2, 2, 1, 1, 2, 2, 2, 2, 1, 1, -1, -2, 1,
-1, 1, 2, 1, -1, 2, -1, 1, 2, 1, -2, 2, -2}\\\hline
{\small 1, 1, 1, 1, 1, 1, 1, 1, 1, 1, 2, 2, 1, 1, 2, 2, 2, 2, 1, 2, -2, -2, 1,
2, 1, -1, 1, 2, -1, -1, 2, 2, -2, -1, 1, 2}\\\hline
{\small 1, 1, 1, 1, 1, 1, 1, 1, 1, 1, 2, 2, 1, 1, 2, 2, 2, 2, 1, 2, -1, -1, 1,
-1, 2, 1, -1, -2, 2, 2, 2, 2, -1, -2, 2, 1}\\\hline
{\small 1, 1, 1, 1, 1, 1, 1, 1, 1, 1, 2, 2, 1, 1, 2, 2, 2, 2, 1, 2, 2, 2, 1,
2, 1, -2, -1, -2, 1, -1, 2, -2, 1, -1, 2, -1}\\\hline
{\small 1, 1, 1, 1, 1, 1, 1, 1, 2, 2, 2, 2, 1, 2, 1, 2, -1, -1, 1, 2, 1, 2,
-2, -2, 2, 2, 1, 2, 1, 2, 2, 2, -1, 1, -2, -1}\\\hline
{\small 1, 1, 1, 1, 1, 2, 1, 1, 1, 1, 2, 1, 1, 1, 2, -1, 2, -2, 1, 1, 2, -2,
1, -2, 1, 2, 2, -1, 1, 1, 1, 2, 2, -2, 1, -1}\\\hline
{\small 1, 1, 1, 1, 1, 2, 1, 1, 1, 2, -2, 1, 1, 2, 2, 2, -1, 2, 1, 2, 2, 2,
-2, 1, 1, -1, -2, 1, 2, 2, 1, -1, -2, 2, 1, 2}\\\hline
{\small 1, 1, 1, 1, 1, 2, 1, 1, 1, 2, -1, 2, 1, 1, 2, -2, 1, -2, 1, 1, 2, -2,
2, -1, 1, 2, 2, 1, -1, 1, 1, 2, 2, 2, -1, 2}\\\hline
{\small 1, 1, 1, 1, 1, 2, 1, 1, 1, 2, 2, 2, 1, 1, 1, -1, -2, 1, 1, 1, -1, 2,
2, 1, 1, 2, -1, -1, -2, 1, 1, 2, -2, 1, 1, 1}\\\hline
{\small 1, 1, 1, 1, 1, 2, 1, 1, 1, 2, 2, 2, 1, 1, 2, -1, -2, 2, 1, -1, 2, 2,
2, 2, 1, -1, 2, -1, -2, 1, 1, -2, 2, 1, 1, 1}\\\hline
{\small 1, 1, 1, 1, 1, 2, 1, 1, 1, 2, 2, 2, 1, 2, -1, 1, 2, 1, 1, 2, -1, 2, 2,
-1, 1, 2, -2, 1, 1, -2, 1, 2, -2, 1, 2, -1}\\\hline
\end{tabular}
\
\]

\noindent All $N$s have probabilities $2/371.$%

\[%
\begin{tabular}
[c]{|l|}\hline
Matrices $N_{112},N_{113},\ldots,N_{131}$ containing exactly $8$ negative
entries (first part):\\\hline
{\small 1, 1, 1, 1, 1, 2, 1, 1, 1, 1, 2, 1, 1, 1, 1, 2, 2, 2, 1, 2, 2, 1, 1,
2, 1, -1, -2, -2, 1, -2, 2, -1, -2, -2, 2, -2}\\\hline
{\small 1, 1, 1, 1, 1, 2, 1, 1, 1, 1, 2, 1, 1, 1, 1, 2, 2, 2, 1, 2, 2, 1, 1,
2, 1, -1, -2, -1, 1, -1, 2, -1, -2, -1, 2, -1}\\\hline
{\small 1, 1, 1, 1, 1, 2, 1, 1, 1, 1, 2, 1, 1, 1, 1, 2, 2, 2, 1, 2, 2, 2, 2,
2, 1, -1, -2, -1, 1, -1, 2, -1, -2, -2, 1, -1}\\\hline
{\small 1, 1, 1, 1, 1, 2, 1, 1, 1, 1, 2, 1, 1, 1, 1, 2, 2, 2, 1, 2, 2, 2, 2,
2, 1, -1, -2, -1, 2, -2, 2, -1, -2, -2, 2, -2}\\\hline
{\small 1, 1, 1, 1, 1, 2, 1, 1, 1, 2, -2, -2, 1, 1, 2, 2, -2, -1, 1, 2, -2, 1,
1, -2, 1, 2, -2, 1, 2, -1, 2, 2, 1, 2, 1, 1}\\\hline
{\small 1, 1, 1, 1, 1, 2, 1, 1, 1, 2, -2, -2, 1, 1, 2, 2, -2, -1, 1, 2, 2, 2,
-1, -1, 1, 2, 2, 2, -2, -2, 2, 2, 1, 2, 1, 1}\\\hline
{\small 1, 1, 1, 1, 1, 2, 1, 1, 1, 2, -2, -2, 1, 1, 2, 2, -1, -2, 1, 2, -1, 2,
-2, 1, 1, 2, -2, 2, -2, 2, 2, 2, 2, 2, 1, 2}\\\hline
{\small 1, 1, 1, 1, 1, 2, 1, 1, 1, 2, -2, -2, 1, 1, 2, 2, -1, -2, 1, 2, 1, 2,
-2, -1, 1, 2, 2, 2, -2, -2, 2, 2, 2, 2, 1, 2}\\\hline
{\small 1, 1, 1, 1, 1, 2, 1, 1, 1, 2, 2, 2, 1, 1, -2, 1, 2, -2, 1, 2, -1, 2,
2, -1, 1, 2, -2, 2, 2, -2, 2, 2, -2, 1, 2, -2}\\\hline
{\small 1, 1, 1, 1, 1, 2, 1, 1, 1, 2, 2, 2, 1, 1, -1, 2, -1, 1, 1, 2, -1, 2,
-1, 2, 1, 2, -2, 2, -2, 2, 2, 2, -1, 2, -2, 2}\\\hline
{\small 1, 1, 1, 1, 1, 2, 1, 1, 1, 2, 2, 2, 1, 1, 1, -1, -2, 1, 1, 1, 2, -2,
-2, 1, 1, 2, 2, -1, -2, 1, 2, 1, 2, -2, -2, 2}\\\hline
{\small 1, 1, 1, 1, 1, 2, 1, 1, 1, 2, 2, 2, 1, 1, 1, -1, -2, 1, 1, 2, 2, -1,
-2, 1, 1, -1, -2, 1, 1, 1, 2, -1, -2, 1, 1, 2}\\\hline
{\small 1, 1, 1, 1, 1, 2, 1, 1, 1, 2, 2, 2, 1, 1, 1, -1, -2, 1, 1, 2, 2, -1,
-2, 1, 1, -1, -2, 2, 2, 1, 2, -1, -2, 2, 2, 2}\\\hline
{\small 1, 1, 1, 1, 1, 2, 1, 1, 1, 2, 2, 2, 1, 1, 1, -1, -2, 1, 1, 2, 2, 1, 1,
2, 1, -1, -2, 1, 1, 1, 2, -1, -2, -1, -2, 1}\\\hline
{\small 1, 1, 1, 1, 1, 2, 1, 1, 1, 2, 2, 2, 1, 1, 1, -1, -2, 1, 1, 2, 2, 2, 2,
2, 1, -1, -2, 1, 2, 2, 2, -1, -2, -2, -2, 2}\\\hline
{\small 1, 1, 1, 1, 1, 2, 1, 1, 1, 2, 2, 2, 1, 1, 2, -2, -2, 1, 1, 2, 2, -2,
-2, 2, 1, -2, 1, 1, 2, -2, 2, -2, 2, 1, 2, -2}\\\hline
{\small 1, 1, 1, 1, 1, 2, 1, 1, 1, 2, 2, 2, 1, 1, 2, -2, -2, 1, 1, 2, 2, 1, 1,
2, 1, -2, -2, 1, 2, 1, 2, -2, -1, -2, -1, 1}\\\hline
{\small 1, 1, 1, 1, 1, 2, 1, 1, 1, 2, 2, 2, 1, 1, 2, -1, -2, 2, 1, 2, -1, 2,
2, -1, 1, 2, -2, 2, 2, -2, 2, 2, 2, -1, -2, 2}\\\hline
{\small 1, 1, 1, 1, 1, 2, 1, 1, 1, 2, 2, 2, 1, 1, 2, -1, -2, 2, 1, 2, 1, -2,
-2, 1, 1, 2, 2, -2, -2, 2, 2, 2, 2, -1, -2, 2}\\\hline
{\small 1, 1, 1, 1, 1, 2, 1, 1, 1, 2, 2, 2, 1, 1, 2, -1, -2, 2, 1, 2, 2, -2,
-2, 2, 1, -1, 1, 1, -2, 1, 2, -1, 2, 2, -2, 2}\\\hline
\end{tabular}
\
\]

\noindent Only $N_{116}$, $N_{125}$, $N_{128}$ have probabilities $1/371$; the
other $N$s have probabilities $2/371.$%

\[%
\begin{tabular}
[c]{|l|}\hline
Matrices $N_{132},N_{133},\ldots,N_{151}$ containing exactly $8$ negative
entries (second part):\\\hline
{\small 1, 1, 1, 1, 1, 2, 1, 1, 1, 2, 2, 2, 1, 1, 2, -1, -2, 2, 1, 2, 2, -1,
-1, 2, 1, -2, 1, 1, -1, 1, 2, -2, 2, 1, -2, 2}\\\hline
{\small 1, 1, 1, 1, 1, 2, 1, 1, 1, 2, 2, 2, 1, 1, 2, -1, -2, 2, 1, 2, 2, 2, 2,
2, 1, -2, -1, -2, -2, 1, 2, -2, -2, 1, 2, 2}\\\hline
{\small 1, 1, 1, 1, 1, 2, 1, 1, 1, 2, 2, 2, 1, 1, 2, -1, -2, 2, 1, 2, 2, 2, 2,
2, 1, -1, -1, -1, -2, 1, 2, -1, -2, 2, 2, 2}\\\hline
{\small 1, 1, 1, 1, 1, 2, 1, 1, 1, 2, 2, 2, 1, 1, 2, -1, -1, 1, 1, 2, 2, -1,
-1, 2, 1, -2, 1, 1, 2, -2, 2, -2, 2, 1, 2, -2}\\\hline
{\small 1, 1, 1, 1, 1, 2, 1, 1, 1, 2, 2, 2, 1, 1, 2, -1, -1, 1, 1, 2, 2, 2, 2,
2, 1, -1, -2, 1, 2, 2, 2, -1, -1, -2, -1, 2}\\\hline
{\small 1, 1, 1, 1, 1, 2, 1, 1, 1, 2, 2, 2, 1, 1, 2, 1, 2, 2, 1, 2, -1, -2,
-2, -1, 1, 2, -2, -2, -2, -2, 2, 2, 2, 1, 2, 2}\\\hline
{\small 1, 1, 1, 1, 1, 2, 1, 1, 1, 2, 2, 2, 1, 1, 2, 1, 2, 2, 1, 2, -1, -1,
-1, -1, 1, 2, -1, -2, -1, -2, 2, 2, 1, 2, 2, 2}\\\hline
{\small 1, 1, 1, 1, 1, 2, 1, 1, 1, 2, 2, 2, 1, 1, 2, 1, 2, 2, 1, 2, 2, 2, 2,
2, 1, -2, -1, -2, -2, 1, 2, -2, -2, -1, -2, 2}\\\hline
{\small 1, 1, 1, 1, 1, 2, 1, 1, 1, 2, 2, 2, 1, 1, 2, 1, 2, 2, 1, 2, 2, 2, 2,
2, 1, -1, -1, -2, -1, 1, 2, -1, -2, -1, -1, 2}\\\hline
{\small 1, 1, 1, 1, 1, 2, 1, 1, 1, 2, 2, 2, 1, 1, 2, 2, 2, 1, 1, 2, 2, 2, 2,
2, 1, -2, -2, -1, -2, 1, 2, -2, -1, -1, -2, 1}\\\hline
{\small 1, 1, 1, 1, 1, 2, 1, 1, 1, 2, 2, 2, 1, 1, 2, 2, 2, 1, 1, 2, 2, 2, 2,
2, 1, -2, 1, -1, -2, -2, 2, -2, 2, -1, -2, -2}\\\hline
{\small 1, 1, 1, 1, 1, 2, 1, 1, 1, 2, 2, 2, 1, 1, 2, 2, 2, 1, 1, 2, 2, 2, 2,
2, 1, -1, -2, -1, -2, 2, 2, -1, -1, -1, -2, 2}\\\hline
{\small 1, 1, 1, 1, 1, 2, 1, 1, 1, 2, 2, 2, 1, 1, 2, 2, 2, 1, 1, 2, 2, 2, 2,
2, 1, -1, 1, -1, -2, -1, 2, -1, 2, -1, -2, -1}\\\hline
{\small 1, 1, 1, 1, 1, 2, 1, 1, 1, 2, 2, 2, 1, 2, 2, 2, 2, 2, 2, 1, 2, 2, 2,
2, 2, -2, 1, -1, -2, -1, 2, -2, 2, -1, -2, -2}\\\hline
{\small 1, 1, 1, 2, 2, 2, 1, 1, 2, 1, 2, 2, 1, 1, 2, -2, 1, -2, 1, 2, 2, -1,
2, -1, 1, 2, 2, -2, 2, -2, 2, 2, 2, -1, 2, -2}\\\hline
{\small 1, 1, 1, 2, 2, 2, 1, 2, 2, 1, -1, -2, 1, 2, 2, 2, 2, 2, 2, 2, 2, 2,
-1, -2, 2, -1, -2, 2, 1, 1, 2, -1, -2, 2, 2, 2}\\\hline
{\small 1, 1, 1, 2, 2, 2, 1, 2, 2, 1, 1, 2, 1, 2, 2, 2, 2, 2, 2, 1, 2, 2, 2,
2, 2, -1, 2, -1, -2, -1, 2, -2, 2, -1, -2, -2}\\\hline
{\small 1, 1, 1, 2, 2, 2, 1, 2, 2, 1, 1, 2, 1, 2, 2, 2, 2, 2, 2, 2, 2, -1, -2,
2, 2, -1, -2, 2, 2, 2, 2, -1, -2, -1, -2, 1}\\\hline
{\small 1, 1, 1, 2, 2, 2, 1, 2, 2, 2, 2, 2, 1, -2, -2, 1, 1, 2, 2, 1, -2, 1,
2, -1, 2, 2, -2, 1, 2, -2, 2, -1, -2, 2, 2, 2}\\\hline
{\small 1, 1, 1, 2, 2, 2, 1, 2, 2, 2, 2, 2, 1, -2, -2, 1, 1, 2, 2, 1, 2, 2,
-1, -1, 2, 1, 2, 2, -2, -2, 2, -2, -2, 2, 1, 2}\\\hline
\end{tabular}
\
\]

\noindent Only $N_{149}$, $N_{150}$ have probabilities $1/371$; the other $N$s
have probabilities $2/371.$%

\[%
\begin{tabular}
[c]{|l|}\hline
Matrices $N_{152},N_{153},\ldots,N_{166}$ containing exactly $9$ negative
entries:\\\hline
{\small 1, 1, 1, 1, 1, 1, 1, 1, 1, -1, -1, -2, 1, 1, 2, 2, 2, 2, 1, 2, 1, -1,
-2, -2, 1, 2, 2, -1, -1, -2, 2, 2, 2, 2, 1, 1}\\\hline
{\small 1, 1, 1, 1, 1, 1, 1, 1, 1, 1, 2, -1, 1, 1, 2, -1, -1, 2, 1, 1, 2, -2,
-2, 2, 1, 2, 2, -1, 1, -1, 1, 2, 2, -2, 1, -2}\\\hline
{\small 1, 1, 1, 1, 1, 1, 1, 1, 1, 1, 2, -1, 1, 1, 2, -1, -1, 2, 1, 1, 2, -2,
-1, 1, 1, 2, 2, -1, 2, -2, 1, 2, 2, -2, 1, -2}\\\hline
{\small 1, 1, 1, 1, 1, 1, 1, 1, 1, 1, 2, 2, 1, 1, 1, -1, -1, -2, 1, 1, 2, -1,
-2, -2, 1, 2, 2, 1, 1, 1, 2, 1, 2, -1, -1, -2}\\\hline
{\small 1, 1, 1, 1, 1, 1, 1, 1, 1, 1, 2, 2, 1, 1, 1, -1, -1, -2, 1, 1, 2, -1,
-1, -1, 1, 2, 2, 1, 2, 2, 2, 1, 2, -1, -1, -2}\\\hline
{\small 1, 1, 1, 1, 1, 1, 1, 1, 1, 1, 2, 2, 1, 1, 2, 2, 2, 2, 1, 1, -1, -2,
-1, -2, 1, 2, -1, -2, 1, 1, 1, -1, 2, 2, -1, -2}\\\hline
{\small 1, 1, 1, 1, 1, 1, 1, 1, 1, 1, 2, 2, 1, 1, 2, 2, 2, 2, 1, 1, -1, -2,
-1, -2, 1, 2, -1, -1, 2, 1, 1, -1, 2, 1, -2, -2}\\\hline
{\small 1, 1, 1, 1, 1, 1, 1, 1, 1, 1, 2, 2, 1, 1, 2, 2, 2, 2, 1, 1, -1, -2, 1,
-1, 1, 2, -1, -1, 1, -1, 1, 2, -1, -2, 1, -2}\\\hline
{\small 1, 1, 1, 1, 1, 1, 1, 1, 1, 2, -1, -2, 1, 1, 2, 2, -1, -1, 1, 2, -1, 1,
1, -2, 1, 2, -1, 1, 2, -1, 2, 2, 1, 2, 1, -1}\\\hline
{\small 1, 1, 1, 1, 1, 2, 1, 1, 1, 2, 2, 2, 1, 1, 2, -2, -2, 1, 1, 2, 2, -1,
-1, 2, 2, -1, 1, 1, 2, 1, 2, -2, 1, -2, -1, -1}\\\hline
{\small 1, 1, 1, 1, 1, 2, 1, 1, 1, 2, 2, 2, 1, 1, 2, -1, -1, 1, 1, 2, 2, 1, 1,
2, 1, -1, -1, 1, 2, 1, 1, -2, -1, -2, -1, -1}\\\hline
{\small 1, 1, 1, 1, 1, 2, 1, 1, 1, 2, 2, 2, 1, 1, 2, 1, 2, 2, 1, -1, 1, -2, 1,
-2, 1, -1, 2, -2, 1, -1, 2, -1, 2, -2, 2, -1}\\\hline
{\small 1, 1, 1, 1, 1, 2, 1, 1, 1, 2, 2, 2, 1, 1, 2, 1, 2, 2, 1, -1, 2, -2, 1,
-1, 1, -2, 1, -1, 1, -2, 2, -2, 2, -1, 2, -1}\\\hline
{\small 1, 1, 1, 1, 1, 2, 1, 1, 1, 2, 2, 2, 1, 2, 2, 1, 1, 2, 1, 2, 2, -1, -2,
1, 1, -1, -2, 2, 2, 1, 1, -1, -2, -1, -2, -1}\\\hline
{\small 1, 1, 1, 1, 2, 2, 1, 1, 1, 2, -1, -2, 1, 1, 2, 2, 1, 1, 1, 2, -2, -1,
1, -2, 1, 2, -2, -1, 2, -1, 2, 2, 1, 2, 1, -1}\\\hline
\end{tabular}
\
\]

\noindent Only $N_{160}$, $N_{162}$, $N_{165}$ have probabilities $1/371$; the
other $N$s have probabilities $2/371.$%
\[%
\begin{tabular}
[c]{|l|}\hline
Matrices $N_{167},N_{168},\ldots,N_{178}$ containing exactly $11$ negative
entries:\\\hline
{\small 1, 1, 1, 1, 1, 1, 1, 1, 1, -1, -1, -2, 1, 2, 2, 1, -2, -2, 1, 2, 2, 2,
-2, -1, 2, 1, 2, -2, 1, -2, 2, 1, 2, -2, 2, -1}\\\hline
{\small 1, 1, 1, 1, 1, 1, 1, 1, 1, -1, -1, -2, 1, 2, 2, 1, -1, -1, 1, 2, 2, 1,
-2, -2, 2, 1, 2, -1, 1, -1, 2, 1, 2, -2, 1, -2}\\\hline
{\small 1, 1, 1, 1, 1, 1, 1, 1, 1, 1, 2, 2, 1, 1, 1, -1, -1, -2, 1, 2, 2, 1,
1, 1, 1, -1, -2, -2, 1, -1, 2, -1, -2, -2, 2, -1}\\\hline
{\small 1, 1, 1, 1, 1, 1, 1, 1, 1, 1, 2, 2, 1, 1, 2, 2, 1, 1, 1, -1, 1, 2, -2,
-2, 1, -2, -1, 1, -1, -2, 2, -2, -1, 2, -1, -2}\\\hline
{\small 1, 1, 1, 1, 1, 1, 1, 1, 1, 1, 2, 2, 1, 1, 2, 2, 2, 2, 1, -1, 1, 2, -1,
-1, 1, -2, -1, 1, -1, -2, 2, -2, -1, 2, -1, -2}\\\hline
{\small 1, 1, 1, 1, 1, 1, 1, 1, 2, 2, 2, 2, 1, 1, -1, -1, -2, -2, 1, 2, 2, -1,
1, -1, 1, 2, 2, -2, 1, -2, 2, 2, -1, 1, -2, -1}\\\hline
{\small 1, 1, 1, 1, 1, 2, 1, 1, 1, 2, 2, 2, 1, 1, -1, 1, 2, -1, 1, 1, -2, -1,
1, -2, 1, 2, -1, -2, 1, -1, 1, 2, -2, -2, 1, -2}\\\hline
{\small 1, 1, 1, 1, 1, 2, 1, 1, 1, 2, 2, 2, 1, 1, -1, 1, 2, -1, 1, 1, -2, -1,
1, -2, 1, 2, -1, -1, 2, -1, 1, 2, -2, -1, 2, -2}\\\hline
{\small 1, 1, 1, 1, 1, 2, 1, 1, 1, 2, 2, 2, 1, 1, 1, -1, -2, 1, 1, 1, -1, -2,
-2, 1, 1, 2, -1, -2, -2, 2, 1, 2, -2, -2, -1, 2}\\\hline
{\small 1, 1, 1, 1, 1, 2, 1, 1, 1, 2, 2, 2, 1, 1, 2, -1, -2, 2, 1, -1, 1, -2,
-2, 1, 1, -1, 2, -2, -2, 2, 1, -2, 1, -2, -1, 1}\\\hline
{\small 1, 1, 1, 1, 1, 2, 1, 1, 1, 2, 2, 2, 1, 2, -1, 2, 2, -1, 1, 2, -1, -1,
-2, 1, 1, 2, -2, 1, 1, -2, 1, 2, -2, -1, -2, -1}\\\hline
{\small 1, 1, 1, 1, 2, 2, 1, 1, 2, 2, 1, 1, 1, 1, 2, 2, 2, 2, 1, -1, 1, 2, -2,
-2, 1, -1, -2, -2, 1, 2, 1, -2, -2, -1, -2, -1}\\\hline
\end{tabular}
\ \ \ \
\]

\noindent Only $N_{167}$, $N_{168}$, $N_{178}$ have probabilities $1/371$; the
other $N$s have probabilities $2/371.$%

\[%
\begin{tabular}
[c]{|l|}\hline
Matrices $N_{179},N_{180},\ldots,N_{193}$ containing exactly $12$ negative
entries:\\\hline
{\small 1, 1, 1, 1, 1, 2, 1, 1, 1, 2, 2, 2, 1, 1, -2, -1, -2, -2, 1, 2, -1, 2,
2, -1, 1, 2, -2, 2, 2, -2, 2, 2, -2, -1, -2, -2}\\\hline
{\small 1, 1, 1, 1, 1, 2, 1, 1, 1, 2, 2, 2, 1, 1, -2, 1, 2, -2, 1, 2, -1, -2,
-2, -1, 1, 2, -2, -2, -2, -2, 2, 2, -2, 1, 2, -2}\\\hline
{\small 1, 1, 1, 1, 1, 2, 1, 1, 1, 2, 2, 2, 1, 1, -1, 1, 2, -1, 1, 2, -1, -2,
-2, -1, 1, 2, -2, -2, -2, -2, 2, 2, -1, 1, 2, -1}\\\hline
{\small 1, 1, 1, 1, 1, 2, 1, 1, 1, 2, 2, 2, 1, 1, -1, 1, 2, -1, 1, 2, -1, -1,
-1, -1, 1, 2, -2, -1, -1, -2, 2, 2, -1, 1, 2, -1}\\\hline
{\small 1, 1, 1, 1, 1, 2, 1, 1, 1, 2, 2, 2, 1, 1, 1, -1, -2, 1, 1, 2, 2, -2,
-2, 2, 1, -1, -2, 1, -2, -2, 2, -1, -2, 2, -2, -2}\\\hline
{\small 1, 1, 1, 1, 1, 2, 1, 1, 1, 2, 2, 2, 1, 1, 1, -1, -2, 1, 1, 2, 2, -1,
-2, 1, 1, -1, -2, -2, -2, 1, 2, -1, -2, -2, -2, 2}\\\hline
{\small 1, 1, 1, 1, 1, 2, 1, 1, 1, 2, 2, 2, 1, 1, 1, -1, -2, 1, 1, 2, 2, -1,
-2, 1, 1, -1, -2, -1, -1, 1, 2, -1, -2, -1, -1, 2}\\\hline
{\small 1, 1, 1, 1, 1, 2, 1, 1, 1, 2, 2, 2, 1, 1, 1, -1, -2, 1, 1, 2, 2, -1,
-1, 2, 1, -1, -2, 1, -1, -1, 2, -1, -2, 1, -2, -1}\\\hline
{\small 1, 1, 1, 1, 1, 2, 1, 1, 1, 2, 2, 2, 1, 1, 2, -2, -2, 1, 1, 2, 2, -2,
-2, 2, 1, -2, 1, -1, -2, -2, 2, -2, 2, -1, -2, -2}\\\hline
{\small 1, 1, 1, 1, 1, 2, 1, 1, 1, 2, 2, 2, 1, 1, 2, -1, -2, 2, 1, 2, -1, -2,
-2, -1, 1, 2, -2, -2, -2, -2, 2, 2, 2, -1, -2, 2}\\\hline
{\small 1, 1, 1, 1, 1, 2, 1, 1, 1, 2, 2, 2, 1, 1, 2, -1, -2, 2, 1, 2, -1, -1,
-1, -1, 1, 2, -2, -1, -1, -2, 2, 2, 2, -1, -2, 2}\\\hline
{\small 1, 1, 1, 1, 1, 2, 1, 1, 1, 2, 2, 2, 1, 1, 2, -1, -1, 1, 1, 2, 2, -1,
-1, 2, 1, -2, 1, -1, -2, -2, 2, -2, 2, -1, -2, -2}\\\hline
{\small 1, 1, 1, 1, 1, 2, 1, 1, 1, 2, 2, 2, 1, 1, 2, -1, -1, 1, 1, 2, 2, -1,
-1, 2, 1, -1, 1, -1, -2, -1, 2, -1, 2, -1, -2, -1}\\\hline
{\small 1, 1, 1, 2, 2, 2, 1, 2, 2, 1, -1, -2, 1, 2, 2, 2, 2, 2, 2, 2, 2, 2,
-1, -2, 2, -1, -2, 2, -1, -1, 2, -1, -2, 2, -2, -2}\\\hline
{\small 1, 1, 1, 2, 2, 2, 1, 2, 2, 2, 2, 2, 1, -2, -2, 1, -1, -2, 2, -1, -2,
2, 1, 1, 2, -1, -2, 2, 2, 2, 2, -2, -2, 2, -1, -2}\\\hline
\end{tabular}
\
\]

\noindent Only $N_{186}$ has probability $1/371$; $N_{182}$, $N_{185}$,
$N_{191}$, $N_{192}$, $N_{193}$ have probabilities $2/1113$; the other $N$s
have probabilities $2/371.$%
\[%
\begin{tabular}
[c]{|l|}\hline
Matrices $N_{194},N_{195},N_{196}$ containing exactly $13$ negative
entries:\\\hline
{\small 1, 1, 1, 1, 1, 1, 1, 1, 2, 2, 2, -1, 1, -2, -1, -1, -2, 1, 2, -1, 1,
2, 1, -1, 2, -2, 1, 2, 1, -2, 2, -2, -1, -1, -2, 2}\\\hline
{\small 1, 1, 1, 1, 1, 2, 1, 1, 1, 2, -1, 2, 1, 1, -1, 2, -2, -1, 1, 1, -2, 1,
-1, -2, 1, 2, -1, 2, -1, -1, 1, 2, -2, 2, -1, -2}\\\hline
{\small 1, 1, 1, 1, 1, 2, 1, 1, 2, 2, -1, 1, 1, -1, 1, 2, -2, -1, 1, -1, -2,
-1, 2, 1, 1, -2, -1, 1, -1, -1, 2, -1, -1, 1, 1, 1}\\\hline
\end{tabular}
\ \
\]

\noindent$\mathbb{P}(\Xi=N_{194})=1/371$, $\mathbb{P}(\Xi=N_{195})=2/371$,
$\mathbb{P}(\Xi=N_{196})=1/1113.$%

\[%
\begin{tabular}
[c]{|l|}\hline
Matrices $N_{197},N_{198},N_{199}$ containing exactly $15$ negative
entries:\\\hline
{\small 1, 1, 1, 1, 1, 2, 1, 1, 1, -1, -2, -2, 1, 1, 2, -1, -2, -1, 1, 2, 1,
-2, -1, -1, 1, 2, 2, -2, -2, -1, 2, 2, 2, -1, -2, -1}\\\hline
{\small 1, 1, 1, 1, 1, 2, 1, 1, 1, -1, -2, -2, 1, 1, 2, -1, -1, -2, 1, 2, 1,
-1, -2, -1, 1, 2, 2, -1, -2, -2, 2, 2, 2, -1, -2, -1}\\\hline
{\small 1, 1, 1, 1, 1, 2, 1, 1, 1, -1, -2, -2, 1, 2, 2, -1, -1, -1, 1, 2, 2,
-1, -2, -2, 1, -1, -2, 1, -1, -1, 1, -1, -2, 2, 1, 2}\\\hline
\end{tabular}
\
\]

\noindent$\mathbb{P}(\Xi=N_{197})=\mathbb{P}(\Xi=N_{198})=2/371$,
$\mathbb{P}(\Xi=N_{199})=1/371.$\pagebreak%

\[%
\begin{tabular}
[c]{|l|}\hline
Matrices $N_{200},N_{201},N_{202},N_{203}$ containing exactly $16$ negative
entries:\\\hline
{\small 1, 1, 1, 1, 1, 2, 1, 1, 1, 2, 2, 2, 1, 1, -2, -1, -2, -2, 1, 2, -1,
-2, -2, -1, 1, 2, -2, -2, -2, -2, 2, 2, -2, -1, -2, -2}\\\hline
{\small 1, 1, 1, 2, 2, 2, 1, 1, 2, 1, -2, -2, 1, -2, 1, -2, -1, -2, 2, -1, 2,
-1, -1, -2, 2, -2, 1, -1, 2, 1, 2, -2, 2, -2, -1, -2}\\\hline
{\small 1, 1, 1, 2, 2, 2, 1, 1, 2, 1, 2, 2, 1, -2, 1, -2, -1, -2, 2, -1, -2,
2, -1, -1, 2, -2, 1, -1, -1, -2, 2, -2, -2, 1, -2, -2}\\\hline
{\small 1, 1, 1, 2, 2, 2, 1, 1, 2, 1, 2, 2, 1, -2, 1, -2, -1, -2, 2, -1, 2,
-2, -1, -1, 2, -2, 1, -2, -2, -2, 2, -2, 2, -2, -1, -2}\\\hline
\end{tabular}
\
\]

\noindent$\mathbb{P}(\Xi=N_{200})=\mathbb{P}(\Xi=N_{203})=2/371$,
$\mathbb{P}(\Xi=N_{201})=\mathbb{P}(\Xi=N_{202})=1/371.\bigskip\bigskip
\bigskip$

\noindent This entire section is based on a conjecture that the list of $203$
canonical representative matrices is complete. \ In other words, we presume
that, for every minimal unimodular zerofree $6\times6$ matrix, there exists an
equivalent matrix in our list. \ The numerical evidence supporting this
conjecture is solid -- a $204^{\text{th}}$ matrix would possess extremely low
probability -- but a rigorous proof that such an event cannot possibly occur
is beyond us.

\section{Addendum III}

Other equivalence classes for minimal unimodular zerofree $7\times7$ matrices include%

\[%
\begin{array}
[c]{ccc}%
N_{3}=\left(
\begin{array}
[c]{ccccccc}%
1 & 1 & 1 & 1 & 1 & 1 & 2\\
1 & 1 & 1 & 1 & 2 & 2 & 2\\
1 & 1 & 1 & 2 & 1 & 2 & 2\\
1 & 1 & 2 & 1 & 1 & -1 & -1\\
1 & 2 & 2 & 2 & 2 & 1 & 1\\
2 & 1 & 2 & 2 & 2 & 1 & 1\\
2 & 2 & 2 & 2 & 2 & 1 & 2
\end{array}
\right)  , &  & N_{3}^{-1}=\left(
\begin{array}
[c]{ccccccc}%
-2 & 1 & 1 & 1 & -2 & -1 & 2\\
-2 & 1 & 1 & 1 & -1 & -2 & 2\\
1 & 1 & 1 & 2 & -1 & -1 & -1\\
1 & -2 & -1 & -2 & 2 & 2 & -1\\
1 & -1 & -2 & -2 & 2 & 2 & -1\\
-2 & 2 & 2 & 2 & -2 & -2 & 1\\
2 & -1 & -1 & -1 & 1 & 1 & -1
\end{array}
\right)  ;
\end{array}
\]%
\[%
\begin{array}
[c]{ccc}%
N_{4}=\left(
\begin{array}
[c]{ccccccc}%
1 & 1 & 1 & 1 & 1 & 1 & 2\\
1 & 1 & 1 & 1 & 1 & 2 & 1\\
1 & 1 & 1 & 1 & 2 & 2 & 2\\
1 & 1 & 1 & 2 & 1 & 2 & 2\\
1 & 1 & 2 & 1 & 1 & 2 & 2\\
1 & 2 & 1 & 1 & 1 & 2 & 2\\
1 & 2 & 2 & -1 & -1 & 1 & 1
\end{array}
\right)  , &  & N_{4}^{-1}=\left(
\begin{array}
[c]{ccccccc}%
1 & 1 & 1 & 1 & -2 & -2 & 1\\
1 & 1 & -2 & -2 & 1 & 2 & -1\\
1 & 1 & -2 & -2 & 2 & 1 & -1\\
1 & 1 & -2 & -1 & 1 & 1 & -1\\
1 & 1 & -1 & -2 & 1 & 1 & -1\\
-2 & -1 & 2 & 2 & -1 & -1 & 1\\
-1 & -2 & 2 & 2 & -1 & -1 & 1
\end{array}
\right)  ;
\end{array}
\]
\bigskip%
\[%
\begin{array}
[c]{ccc}%
N_{5}=N_{3}^{\prime}; &  & N_{6}=N_{4}^{\prime};
\end{array}
\]
\
\[%
\begin{array}
[c]{ccc}%
N_{7}=\left(
\begin{array}
[c]{ccccccc}%
1 & 1 & 1 & 1 & 1 & 1 & 2\\
1 & 1 & 1 & 1 & 1 & 2 & 1\\
1 & 1 & 1 & 1 & 2 & 2 & 2\\
1 & 1 & 1 & 2 & 1 & 2 & 2\\
1 & 1 & 2 & 1 & 1 & 2 & 2\\
1 & 2 & 1 & 1 & -2 & -1 & -1\\
1 & 2 & 2 & 2 & -1 & 1 & 1
\end{array}
\right)  , &  & N_{7}^{-1}=\left(
\begin{array}
[c]{ccccccc}%
1 & 1 & -2 & 1 & 1 & 1 & -2\\
1 & 1 & 1 & -2 & -2 & -1 & 2\\
-2 & -2 & 1 & 1 & 2 & 1 & -1\\
-2 & -2 & 1 & 2 & 1 & 1 & -1\\
-2 & -2 & 2 & 1 & 1 & 1 & -1\\
1 & 2 & -1 & -1 & -1 & -1 & 1\\
2 & 1 & -1 & -1 & -1 & -1 & 1
\end{array}
\right)  ;
\end{array}
\]
\bigskip\bigskip%
\[%
\begin{array}
[c]{ccc}%
N_{8}=\left(
\begin{array}
[c]{ccccccc}%
1 & 1 & 1 & 1 & 1 & 2 & 2\\
1 & 1 & 1 & 1 & 2 & 1 & 2\\
1 & 1 & 1 & 2 & 2 & 2 & 2\\
1 & 1 & 1 & -2 & 2 & -1 & 1\\
1 & 1 & 2 & 1 & 2 & 2 & 2\\
1 & 2 & 1 & 1 & 2 & 2 & 2\\
2 & 1 & 1 & -2 & 2 & -1 & 2
\end{array}
\right)  , &  & N_{8}^{-1}=\left(
\begin{array}
[c]{ccccccc}%
-2 & -2 & 1 & -2 & 1 & 1 & 2\\
-2 & 1 & -2 & -2 & 1 & 2 & 1\\
-2 & 1 & -2 & -2 & 2 & 1 & 1\\
-2 & 1 & -1 & -2 & 1 & 1 & 1\\
1 & -1 & 2 & 2 & -1 & -1 & -1\\
2 & -2 & 2 & 2 & -1 & -1 & -1\\
2 & 2 & -1 & 1 & -1 & -1 & -1
\end{array}
\right)  ;
\end{array}
\]
\bigskip\bigskip%
\[%
\begin{array}
[c]{ccc}%
N_{9}=\left(
\begin{array}
[c]{ccccccc}%
1 & 1 & 1 & 1 & 1 & 1 & 2\\
1 & 1 & 1 & 1 & 1 & 2 & 1\\
1 & 1 & 1 & 1 & 2 & 1 & 1\\
1 & 1 & 2 & -2 & 1 & 1 & -2\\
1 & 2 & 2 & 2 & 2 & 2 & 2\\
2 & 1 & 2 & -1 & 2 & 2 & -1\\
2 & 2 & 2 & 1 & 2 & 2 & 2
\end{array}
\right)  , &  & N_{9}^{-1}=\left(
\begin{array}
[c]{ccccccc}%
-2 & -2 & -2 & -2 & 1 & 2 & 2\\
-2 & 1 & 1 & 1 & -1 & -2 & 2\\
1 & -2 & -2 & -1 & 2 & 2 & -1\\
-2 & -2 & -2 & -2 & 2 & 2 & 1\\
1 & 1 & 2 & 1 & -1 & -1 & -1\\
1 & 2 & 1 & 1 & -1 & -1 & -1\\
2 & 1 & 1 & 1 & -1 & -1 & -1
\end{array}
\right)  ;
\end{array}
\]
\bigskip\bigskip%
\[%
\begin{array}
[c]{ccccc}%
N_{10}=N_{7}^{\prime}; &  & N_{11}=N_{8}^{\prime}; &  & N_{12}=N_{9}^{\prime};
\end{array}
\]
\bigskip\pagebreak%
\[%
\begin{array}
[c]{ccc}%
N_{13}=\left(
\begin{array}
[c]{ccccccc}%
1 & 1 & 1 & 1 & 1 & 2 & 2\\
1 & 1 & 1 & 1 & 2 & 1 & 2\\
1 & 1 & 1 & 2 & 2 & 2 & 2\\
1 & 1 & 2 & 1 & 2 & 2 & 2\\
1 & 2 & 1 & 1 & 2 & 2 & 2\\
1 & 2 & -1 & -1 & 1 & -2 & 1\\
2 & 2 & -1 & -1 & 1 & -2 & 2
\end{array}
\right)  , &  & N_{13}^{-1}=\left(
\begin{array}
[c]{ccccccc}%
-2 & -2 & 1 & 1 & 1 & -2 & 2\\
1 & -2 & 1 & 1 & -1 & 2 & -1\\
1 & -2 & 1 & 2 & -2 & 2 & -1\\
1 & -2 & 2 & 1 & -2 & 2 & -1\\
-2 & 2 & -1 & -1 & 2 & -2 & 1\\
-1 & 1 & -1 & -1 & 2 & -2 & 1\\
2 & 2 & -1 & -1 & -1 & 1 & -1
\end{array}
\right)  ;
\end{array}
\]
\bigskip%
\[%
\begin{array}
[c]{ccc}%
N_{14}=\left(
\begin{array}
[c]{ccccccc}%
1 & 1 & 1 & 1 & 1 & 2 & 2\\
1 & 1 & 1 & 2 & -2 & 1 & -1\\
1 & 2 & 2 & 2 & -1 & 2 & 1\\
2 & 1 & 2 & 2 & -1 & 2 & 1\\
2 & 2 & 1 & 2 & -1 & 2 & 1\\
2 & 2 & 2 & 2 & -1 & 2 & 2\\
2 & 2 & 2 & 2 & -2 & 2 & 1
\end{array}
\right)  , &  & N_{14}^{-1}=\left(
\begin{array}
[c]{ccccccc}%
-2 & -2 & 1 & 2 & 2 & -1 & -1\\
-2 & -2 & 2 & 1 & 2 & -1 & -1\\
-2 & -2 & 2 & 2 & 1 & -1 & -1\\
1 & 2 & -1 & -1 & -1 & 2 & -1\\
-2 & -2 & 2 & 2 & 2 & -1 & -2\\
2 & 1 & -1 & -1 & -1 & -1 & 2\\
2 & 2 & -2 & -2 & -2 & 2 & 1
\end{array}
\right)  ;
\end{array}
\]
\bigskip%
\[%
\begin{array}
[c]{ccc}%
N_{15}=\left(
\begin{array}
[c]{ccccccc}%
1 & 1 & 1 & 1 & 1 & 1 & 2\\
1 & 1 & 1 & 1 & 2 & 2 & 2\\
1 & 1 & 2 & -2 & 1 & -1 & -1\\
1 & 2 & 2 & -1 & 2 & 1 & 1\\
2 & 1 & 2 & -1 & 2 & 1 & 1\\
2 & 2 & 2 & -1 & 2 & 1 & 2\\
2 & 2 & 2 & -2 & 2 & 1 & 1
\end{array}
\right)  , &  & N_{15}^{-1}=\left(
\begin{array}
[c]{ccccccc}%
1 & 1 & 1 & -2 & -1 & -1 & 2\\
1 & 1 & 1 & -1 & -2 & -1 & 2\\
1 & -2 & -1 & 2 & 2 & -1 & -1\\
1 & 1 & 1 & -1 & -1 & -1 & 1\\
-2 & 2 & 1 & -1 & -1 & 2 & -1\\
1 & -1 & -1 & 1 & 1 & -2 & 1\\
-1 & -1 & -1 & 1 & 1 & 2 & -2
\end{array}
\right)  ;
\end{array}
\]
\bigskip%
\[%
\begin{array}
[c]{ccc}%
N_{16}=\left(
\begin{array}
[c]{ccccccc}%
1 & 1 & 1 & 1 & 1 & 1 & 2\\
1 & 1 & 1 & 1 & 1 & 2 & 1\\
1 & 1 & 1 & -1 & -2 & 1 & 1\\
1 & 1 & 2 & -2 & -2 & 1 & 1\\
1 & 2 & 2 & 2 & 2 & 2 & 2\\
2 & 1 & 2 & -1 & -1 & 2 & 2\\
2 & 2 & 2 & -1 & -2 & 2 & 2
\end{array}
\right)  , &  & N_{16}^{-1}=\left(
\begin{array}
[c]{ccccccc}%
-2 & -2 & -2 & -2 & 1 & 2 & 2\\
1 & 1 & -2 & 1 & -1 & -2 & 2\\
-2 & -2 & 1 & -1 & 2 & 2 & -1\\
-2 & -2 & 2 & -2 & 2 & 2 & -1\\
1 & 1 & -2 & 1 & -1 & -1 & 1\\
1 & 2 & 1 & 1 & -1 & -1 & -1\\
2 & 1 & 1 & 1 & -1 & -1 & -1
\end{array}
\right)  ;
\end{array}
\]
\bigskip%
\[%
\begin{array}
[c]{ccc}%
N_{17}=\left(
\begin{array}
[c]{ccccccc}%
1 & 1 & 1 & 1 & 1 & 1 & 2\\
1 & 1 & 1 & 1 & 2 & 2 & 2\\
1 & 1 & 1 & 2 & 1 & 2 & 2\\
1 & 2 & -1 & 2 & -1 & 1 & 1\\
2 & 1 & -1 & 2 & -1 & 1 & 1\\
2 & 2 & -1 & 2 & -1 & 1 & 2\\
2 & 2 & -2 & 2 & -1 & 1 & 1
\end{array}
\right)  , &  & N_{17}^{-1}=\left(
\begin{array}
[c]{ccccccc}%
1 & 1 & -2 & 1 & 2 & -1 & -1\\
1 & 1 & -2 & 2 & 1 & -1 & -1\\
1 & 1 & -2 & 2 & 2 & -1 & -2\\
1 & -2 & 2 & -1 & -1 & -1 & 2\\
1 & -1 & 1 & -1 & -1 & -1 & 2\\
-2 & 2 & -1 & 1 & 1 & 1 & -2\\
-1 & -1 & 2 & -2 & -2 & 2 & 1
\end{array}
\right)  ;
\end{array}
\]%
\[%
\begin{array}
[c]{ccc}%
N_{18}=\left(
\begin{array}
[c]{ccccccc}%
1 & 1 & 1 & 1 & 1 & 1 & 2\\
1 & 1 & 1 & 1 & 1 & 2 & 1\\
1 & 1 & 1 & 2 & -2 & 1 & -2\\
1 & 1 & 2 & 1 & -2 & 1 & -2\\
1 & 2 & 2 & 2 & -1 & 2 & -1\\
2 & 1 & 2 & 2 & -1 & 2 & -1\\
2 & 2 & 2 & 2 & 1 & 2 & 2
\end{array}
\right)  , &  & N_{18}^{-1}=\left(
\begin{array}
[c]{ccccccc}%
-2 & 1 & 1 & 1 & -2 & -1 & 2\\
-2 & 1 & 1 & 1 & -1 & -2 & 2\\
1 & -2 & -2 & -1 & 2 & 2 & -1\\
1 & -2 & -1 & -2 & 2 & 2 & -1\\
-2 & -2 & -2 & -2 & 2 & 2 & 1\\
1 & 2 & 1 & 1 & -1 & -1 & -1\\
2 & 1 & 1 & 1 & -1 & -1 & -1
\end{array}
\right)  ;
\end{array}
\]
\bigskip%
\[
N_{19}=N_{18}^{\prime};
\]
\bigskip%
\[%
\begin{array}
[c]{ccc}%
N_{20}=\left(
\begin{array}
[c]{ccccccc}%
1 & 1 & 1 & 1 & 1 & 1 & 2\\
1 & 1 & 1 & 1 & 1 & 2 & 1\\
1 & 1 & 1 & -1 & -2 & 1 & -2\\
1 & 1 & 2 & -2 & -2 & 1 & -2\\
1 & 2 & 2 & -1 & -1 & 2 & -1\\
2 & 1 & 2 & -1 & -1 & 2 & -1\\
2 & 2 & 2 & -1 & -2 & 2 & -1
\end{array}
\right)  , &  & N_{20}^{-1}=\left(
\begin{array}
[c]{ccccccc}%
1 & -2 & 1 & -2 & 1 & 2 & -1\\
1 & -2 & 1 & -2 & 2 & 1 & -1\\
1 & 1 & 1 & 2 & -1 & -1 & -1\\
1 & 1 & 2 & 1 & -1 & -1 & -1\\
1 & -2 & 1 & -2 & 2 & 2 & -2\\
-2 & 2 & -2 & 1 & -1 & -1 & 2\\
-1 & 1 & -2 & 1 & -1 & -1 & 2
\end{array}
\right)  .
\end{array}
\]

\bigskip

\noindent We believe more matrices await discovery! \ Transposition does not
necessarily create a new equivalence class, e.g., $N_{14}^{\prime}\sim N_{14}%
$, $N_{15}^{\prime}\sim N_{15}$, $N_{16}^{\prime}\sim N_{13}$ and
$N_{20}^{\prime}\sim N_{17}$. \ As a final observation, the matrix%

\[
\tilde{N}=\left(
\begin{array}
[c]{ccccccc}%
1 & 1 & 1 & 1 & 2 & 2 & -1\\
1 & 1 & 1 & 2 & 1 & -1 & 2\\
1 & 2 & 2 & 2 & 2 & 1 & 1\\
2 & 1 & 2 & 2 & 2 & 1 & 1\\
2 & 2 & 1 & 2 & 2 & 1 & 1\\
2 & 2 & 2 & 2 & 2 & 1 & 2\\
2 & 2 & 2 & 2 & 2 & 2 & 1
\end{array}
\right)
\]
\medskip

\noindent bears a certain resemblance to $N_{3}$ and $N_{4}$: all three
matrices possess exactly two negative entries. \ For $\tilde{N}$, however, the
two entries occupy different rows and different columns. \ It can be shown
that $\tilde{N}=\operatorname*{CanonicalizeMatrix}_{+}[\tilde{N}]$ in the
event we restrict $P$ \&\ $Q$ to the subclass of unsigned permutation
matrices. \ But if we remove the restriction and allow $P$ \&\ $Q$ to be
\textit{signed} permutation matrices, then it follows that
$\operatorname*{CanonicalizeMatrix}[\tilde{N}]=N_{14}$. \ The fact that
$N_{14}$ has no visible similarity with $N_{3}$ or $N_{4}$ is unsatisfactory.
\ The class of signed permutation matrices is perhaps too large and our
definition of canonical representative is maybe too abstruse. \ Some thought
is needed -- resolving such issues -- as a prelude to further work in
unimodular zerofree matrix classification.
\end{document}